\documentclass[11pt,reqno]{amsart} 
\usepackage{amscd,amsmath,amssymb,euscript,graphics}
\usepackage[frame,cmtip,curve,arrow,matrix,line,graph]{xy}

\usepackage[usenames,dvipsnames]{color}
\usepackage{hyperref}
\hypersetup{
colorlinks=true,       				
linkcolor=blue,			   		
citecolor=Magenta}        			

\title{Stokes factors and multilogarithms}
\author{Tom Bridgeland}
\address{
All Souls College,
Oxford,
UK}
\email{bridgeland@maths.ox.ac.uk}
\thanks{T.B. supported by a Royal Society University Research Fellowship}
\author{Valerio Toledano Laredo}
\address{Department of Mathematics,
Northeastern University,
567 Lake Hall,
360 Huntington Avenue,
Boston MA 02115.}
\email{V.ToledanoLaredo@neu.edu}
\thanks{V.T.L. supported in part by NSF grants DMS--0707212 and DMS--0635607}

\newtheorem*{thm}{Theorem}

\newtheorem*{cor}{Corollary}
\newtheorem*{prop}{Proposition}
\newtheorem*{lemma}{Lemma}

\newenvironment{pf}{\paragraph{{\sc Proof}}}{\qed\par\medskip}
\theoremstyle{definition} 
\newtheorem*{defn}{Definition}
\newtheorem*{remark}{Remark}

\newcommand{\GL}{\operatorname{GL}}
\renewcommand{\leq}{\leqslant}
\renewcommand{\geq}{\geqslant}
\newcommand{\ad}{\operatorname{ad}}

\newcommand{\PT}{\operatorname{PT}}

\newcommand{\h}{\mathfrak{h}}
\renewcommand{\t}{\mathfrak t}
\newcommand{\End}{\operatorname{End}}
\newcommand{\D}{\operatorname{\mathcal{D}}}

\newcommand{\nablah}{\widehat{\nabla}}

\newcommand{\C}{\mathbb C}

\newcommand{\blob}{{\scriptscriptstyle\bullet}}
\newcommand{\A}{\mathcal A}

\newcommand{\into}{\hookrightarrow}
\newcommand{\id}{\operatorname{id}}

\newcommand{\Hom}{\operatorname{Hom}}

\newcommand{\tensor}{\otimes}
\newcommand{\R}{\mathbb{R}}

\newcommand{\reg}{{\operatorname{reg}}}
\renewcommand{\b}{\mathfrak{b}}

\newcommand{\g}{\mathfrak{g}}

\newcommand{\hreg}{\h_\reg}

\renewcommand{\Re}{\operatorname{Re}}

\newcommand{\V}{\mathcal{V}}

\newcommand{\god}{\g_{\text{od}}}
\newcommand{\poles}{\mathcal{P}}

\newcommand {\diag}{\operatorname{diag}}

\newcommand{\gl}{\mathfrak{gl}}

\newcommand{\n}{\mathfrak{n}}
\newcommand{\wh}[1]{\widehat{#1}}

\newcommand{\IC}{\mathbb C}
\newcommand{\IH}{\mathbb H}
\newcommand{\Hr}{i\IH_r}
\newcommand{\IN}{\mathbb N}
\newcommand{\IP}{\mathbb P}

\newcommand{\IR}{\mathbb R}

\newcommand{\IZ}{\mathbb Z}

\renewcommand {\L}{\mathcal L}

\newcommand {\calS}{\mathcal S}
\newcommand {\U}{\mathcal U}
\newcommand {\Z}{\mathcal Z}

\newcommand{\Ker}{\operatorname{Ker}}

\newcommand{\eg}{{\it e.g. }}

\newcommand{\ol}[1]{\overline{#1}}
\newcommand {\ul}[1]{\underline{#1}}
\newcommand{\rhs}{right--hand side }
\newcommand{\lhs}{left--hand side }

\newcommand{\wrt}{with respect to }
\newcommand{\aand}{\qquad\text{and}\qquad}

\newcommand {\ray}[1]{\IR_{>0}\cdot #1}
\newcommand {\Ad}{\operatorname{Ad}}
\newcommand {\SL}{\mathcal L}
\newcommand {\nablahv}{\nablah^V}

\newcommand {\Log}{\operatorname{Log}}
\newcommand {\bfX}{\mathbb{X}}
\newcommand {\rep}{representation }
\newcommand {\clockwise}{\stackrel{\curvearrowright}{\prod}}

\newcommand {\PhiZ}{\Phi^Z}
\renewcommand {\Vec}{\operatorname{Vec}}
\newcommand {\RRep}{\operatorname{Rep}}

\newcommand {\LambdaZ}{\Lambda^Z}

\newcommand {\nablaF}{\wh{\nabla}}

\newcommand {\Sym}{\operatorname{Sym}}

\newcommand {\Cstar}{\bigsqcup_{n\geq 1}(\IC^*)^n}

\newcommand {\NC}{\mathcal{NC}}
\newcommand {\CC}{\mathcal{C}}

\newcommand {\wt}[1]{\widetilde{#1}}

\newcommand {\nbd}{neighboorhood }
\newcommand{\fd}{finite--dimensional }
\newcommand{\RH}{Riemann--Hilbert }

\newcommand{\FL}{Fourier--Laplace }

\newcommand{\BJL}{Balser--Jurkat--Lutz }
\newcommand{\BCH}{Baker--Campbell--Hausdorff }
\newcommand {\respnil}{(resp. nilpotent) }
\newcommand {\respuni}{(resp. unipotent) }

\begin{document}

\begin{abstract}
Let $G$ be a complex, affine algebraic group and $\nabla$
a meromorphic connection on the trivial $G$--bundle over
$\IP^1$, with a pole of order 2 at zero and a pole of order
1 at infinity. We show that the map $\calS$ taking the residue
of $\nabla$ at zero to the corresponding Stokes factors is
given by an explicit, universal Lie series whose coefficients
are multilogarithms. Using a non--commutative analogue
of the compositional inversion of formal power series, we
show that the same holds for the inverse of $\calS$, and that
the corresponding Lie series coincides with the generating
function for counting invariants in abelian categories constructed
by D. Joyce.
\end{abstract}
\maketitle

\setcounter{tocdepth}{1}
\tableofcontents

\section{Introduction}

\subsection{}\label{ss:intro intro} 

Let $\V$ be the trivial, rank $n$ complex vector bundle on $\IP^1$
and consider a meromorphic connection on $\V$ of the form
\begin{equation}\label{eq:reference}
\nabla=
d-\left(\frac{Z}{t^2}+\frac{f}{t}\right)dt
\end{equation}
where $Z,f$ are $n\times n$ matrices and $Z$ is diagonal with distinct
eigenvalues $z_1,\ldots,z_n$.

Since $\nabla$ has an irregular singularity at $t=0$, its gauge equivalence
class, as a connection on the unit disk, is determined by its {\it Stokes factors}
\cite{BJL}. These encode the change in the asymptotics of the canonical
fundamental solutions of $\nabla$ across suitable rays in the $t$--plane.
For each such {\it Stokes ray}, that is a ray of the form $\ell=\IR_{>0}(z_i-z_j)$,
the Stokes factor $S_\ell$ is a unipotent matrix whose only non--zero,
off--diagonal entries are of the form $(S_\ell)_{jk}$ with $j,k$ such that
$z_j-z_k\in\ell$.

\subsection{}

The computation of the Stokes factors of $\nabla$ was reduced by \BJL
to that of the analytic continuation of solutions of the \FL transform $\wh
{\nabla}$ of $\nabla$ \cite{BJL2}, that is the Fuchsian connection with
poles at the points $z_1,\ldots,z_n$, given by
\begin{equation}\label{eq:FL ref}
\wh{\nabla}=
d-\sum_{i=1}^n\frac{P_i f}{z-z_i}dz
\end{equation}
where $P_i$ is the projection onto the $Z$--eigenspace corresponding
to $z_i$.

Using the well--known fact that the monodromy of such connections can
be expressed in terms of multilogarithms\footnote{this fact is only implicit
in \cite{BJL2}, as is the above formula for $S_\ell$.} leads in particular to
the following formula for $S_\ell$, with $\ell=\IR_{>0}(z_i-z_j)$ (see
Theorem \ref{one2})
\[S_\ell=
1+\sum_{n\geq 1}\sum_{\substack{1\leq i_1\neq\cdots\neq i_{n+1}\leq n\\z_{i_1}-z_{i_{n+1}}\in\ell}}
M_n(z_{i_1}-z_{i_2},\ldots,z_{i_n}-z_{i_{n+1}})f_{i_1\,i_2}\cdots f_{i_n\,i_{n+1}}\]
where $f_{ij}=P_ifP_j$ is the component of $f$ along the elementary
matrix $E_{ij}$, $f=\sum_{i\neq j}f_{ij}$ is assumed to have zero diagonal
entries, and the function $M_n$ is given by the iterated integral
\[M_n(w_1,\ldots,w_n)=
2\pi i \big.
\int_{[0,w_1+\cdots+w_n]}\frac{dt}{t-w_1}
\circ \cdots \circ \frac{dt}{t-(w_1+\cdots+w_{n-1})}\]

\subsection{}	

One of the goals of this paper is to extend the results of \cite{BJL,BJL2}
in two distinct directions by allowing:
\begin{enumerate}
\item the structure group $GL_n(\IC)$ of the connection $\nabla$ 
to be an arbitrary complex, affine algebraic group $G$.
\item The element $Z$ to be an arbitrary semisimple element of
the Lie algebra of $G$, in particular a diagonal matrix with repeated
eigenvalues for $G=GL_n(\IC)$.
\end{enumerate}
As we explain in \ref{ss:Lie series} and \ref{ss:RH} below, this leads
to new results even when $G=GL_n (\IC)$ and $Z$ has distinct
eigenvalues. In particular, we show that the logarithms of the Stokes
factors are given by {\it universal Lie series}, akin to the \BCH formula,
in the variables $\{f_{ij}\}$. We also obtain an explicit solution of the
corresponding \RH problem, namely the construction of an appropriate
coefficient matrix $f$ starting from prescribed Stokes data, in terms
of such Lie series.

\subsection{}	

With regard to (i), an extension of the \BJL theory of invariants for
meromorphic connections \cite{BJL}, in particular the construction
of canonical fundamental solutions and definition of Stokes data,
was carried out by P. Boalch for an arbitrary complex, reductive
group \cite{Boalch3,Boalch2}. Just as \cite{BJL}, Boalch's treatment
applies to a more general context than the one presented above:
it is local, in that it encompasses meromorphic connections on the
unit disk rather than $\IP^1$, and allows for arbitrary order poles.
It does not however address the extension (ii).

It turns out that an extension of \cite{BJL} to the larger class of affine
algebraic groups may be readily obtained from the corresponding
results for the group $GL_n(\IC)$ by using Tannaka duality, that is
the reconstruction of such groups from their \fd representations. This
requires taking care of the extension (ii) for $G=GL_n(\IC)$ however,
since the assumption that $Z$ be regular semisimple is not stable
under passage to a representation.

When $G=GL_n(\IC)$ and $Z$ is a matrix with repeated eigenvalues,
the extension of \cite{BJL} and \cite{BJL2} was in fact carried out for
arbitrary order poles by \BJL in \cite{BJL3} and \cite{BJL4} respectively\footnote
{We are grateful to Phil Boalch for pointing us to \cite {BJL4}.}. The
main result of \cite {BJL4}, namely the computation of the Stokes
factors in terms of the analytic continuation of suitable associated
functions is not sufficiently explicit for our purposes however. We
therefore give a streamlined treatment which closely parallels that
of \cite{BJL2} and relies on the special form of the connection
$\nabla$.

\subsection{}\label{ss:Lie series} 

Once the canonical fundamental solutions of $\nabla$ are constructed,
one finds, as in \cite{BJL,Boalch3}, that each Stokes factor $S_\ell$ is
a unipotent element of $G$ of the form $S_\ell=\exp(\epsilon_\ell)$,
where $\epsilon_\ell$ lies in the span of the $\ad(Z)$--eigenspaces
corresponding to eigenvalues lying on $\ell$. As $\ell$ varies through
the Stokes rays, the logarithms of the Stokes factors may be assembled
to give a {\it Stokes map} $\calS$ mapping $f$ to $\epsilon=\sum_\ell
\epsilon_\ell$.

Extending the results of \cite{BJL2}, we compute this map explicitly,
first for $G=GL_n(\IC)$ by using the \FL transform and then for an
arbitrary algebraic group by using Tannaka duality. Our results differ
in their form from those of \cite{BJL2} in that by considering the
logarithms of the Stokes factors rather than the Stokes factors
themselves, we express the answer as a universal Lie series
in the components $\{f_\alpha\}$ of $f$ along a root space decomposition
of the Lie algebra of $G$. The coefficients of this series are multilogarithms
evaluated at the eigenvalues of $\ad(Z)$.

\subsection{}\label{ss:RH} 

We then solve the underlying \RH problem, that is the construction of
a connection of the form \eqref{eq:reference} with prescribed Stokes
factors $\{S_\ell\}$, provided these factors are small enough. We do
so by explicitly inverting the Lie series computing the Stokes map
$\calS$, to yield the Taylor series of the local inverse of $\calS$ at
$\epsilon=0$. This inverse is again expressed as an explicit, universal
Lie series involving multilogarithms. To the best of our knowledge,
this result is new even when $G=GL_n(\IC)$ and $Z$ has distinct
eigenvalues.

The inversion of $\calS$ is obtained by using a non--commutative
analogue of the compositional inversion of a formal power series.
The latter expresses the answer as a sum over plane rooted trees
and may be of independent interest.

\subsection{}

Remarkably, the Taylor series of $\calS^{-1}$ coincides with the
generating series for counting invariants in an abelian category
$\A$ constructed by D. Joyce \cite{Joyce}. This allows us to reinterpret
Joyce's construction as the statement that a stability condition on
$\A$ defines Stokes data for a connection of the form \eqref{eq:reference}
with values in the Ringel--Hall Lie algebra of $\A$ \cite{BTL}. Understanding
Joyce's generating series was in fact the main motivation behind
this project, and the initial reason for the need to consider affine
algebraic groups, rather than just $GL_n(\IC)$. Indeed, the groups
underlying such Ringel--Hall Lie algebras are (pro--)solvable.

\subsection{}	

We conclude with a more detailed description of the contents of this paper.
In Section \ref{se:Irregular connections}, we review the definition of the
canonical fundamental solutions of an irregular connection of the form
\eqref{eq:reference} and of the corresponding Stokes data. In Section
\ref{se:imd}, we discuss isomonodromic deformations of such connections.
In Section \ref{se:Stokes map}, we state our results concerning the computation
of the Stokes map and of the Taylor series of its inverse in terms of multilogarithms.
Section \ref{se:multipliers} contains similar results for the Stokes multipliers
of the connection. The rest of the paper contains the proofs of the results of
Section \ref{se:Stokes map}. Specifically, in Section \ref{se:Phil lemma} we
prove the uniqueness of canonical fundamental solutions of the connection
$\nabla$. Section \ref{polylog} covers mostly well--known background material
on the computation of regularised parallel transport for Fuchsian connections
on $\IP^1$ in terms of iterated integrals. In Section \ref{se:Fourier-Laplace},
we prove the existence of the canonical fundamental solutions for the group
$GL_n(\IC)$ and then compute, in Section \ref{se:Stokes factors} the corresponding
Stokes factors in terms of multilogarithms. As in \cite{BJL2}, these results are
derived from those of Section \ref{polylog} by  the use of the \FL transform. In
Section \ref{se:Tannaka}, we prove our main results for an arbitrary affine algebraic
group by using Tannaka duality. Section \ref{comb} gives the non--commutative
generalisation of the compositional inversion of formal power series which is
required to invert the Stokes map.


\section{Irregular connections and Stokes phenomena}
\label{se:Irregular connections}

We review in this section the definition of Stokes data for irregular
connections on $\IP^1$ of the form \eqref{eq:reference}. Our
exposition follows \cite{Boalch3,Boalch2}, where we learnt much
of this material. As pointed out in the Introduction however, we
depart from \cite{Boalch3,Boalch2} and the earlier treatement \cite
{BJL} in that we consider connections having as structure group
an arbitrary complex algebraic group rather than a complex
reductive group.

\subsection{Recollections on algebraic groups} 

We summarise in this paragraph some standard terminology
and facts about algebraic groups and refer the reader to \cite
{Hu} for more details.

By an {\it algebraic group}, we shall always mean an affine
algebraic group $G$ over $\IC$. By a {\it \fd representation}
of $G$, we shall mean a rational representation, that is a
morphism $G\to GL(V)$, where $V$ is a \fd complex vector
space. An algebraic group always possesses a faithful \fd
representation and may therefore be regarded as a {\it linear}
algebraic group, that is a (Zariski) closed subgroup of some
$GL(V)$.

An element $g\in G$ is {\it semisimple} (resp. {\it unipotent})
if it acts by a semisimple \respuni endomorphism on any \fd
representation of $G$. Equivalently, $g$ is semisimple \respuni
if, after embedding $G$ as a closed subgroup of some $GL(V)$,
$g$ is a semisimple \respuni endomorphism of $V$. If $G$ is
semisimple, then $g$ is semisimple \respuni if, and only if,
$\Ad(g)$ is a semisimple \respuni endomorphism of the Lie
algebra $\g$ of $G$.

Similarly, an element $Z\in\g$ is semisimple \respnil if $Z$
acts as a semisimple \respnil endomorphism on any \fd
representation of $G$ or, equivalently, on a faithful
representation of $G$ or, when $G$ is semisimple,
on the adjoint representation of $G$.

\subsection{The irregular connection $\mathbf{\nabla}$}
\label{ss:assumptions}

Let $P$ be the holomorphically trivial, principal $G$--bundle on
$\IP^1$. We shall be concerned with meromorphic connections
on $P$ of the form
\begin{equation}
\label{nab}
\nabla=d-\bigg(\frac{Z}{t^2}+\frac{f}{t}\bigg)dt.
\end{equation}
where $Z,f\in\g$.

Throughout this paper, we assume that the elements $Z,f$ satisfy
the following:
\begin{enumerate}
\item[$(Z)$] $Z$ is semisimple. In particular, $\g$ splits as the
direct sum
\begin{equation}\label{eq:Z decomposition}
\g=\g^Z\oplus[Z,\g]
\end{equation}
where $\g^Z$ is the centraliser of $Z$ and $[Z,\g]$ is the span
of the non--zero eigenspaces of $\ad(Z)$.
\item[$(f)$] The projection of $f$ onto $\g^Z$ corresponding
to the decomposition \eqref{eq:Z decomposition} is zero.
\end{enumerate}
If $G=\GL_n(\C)$ and $Z$ is a diagonal matrix with distinct
eigenvalues, condition ($f$) is the requirement that the diagonal
entries of the matrix $f$ be zero.

We will denote the $\ad(Z)$--eigenspace corresponding to the eigenvalue
$\zeta\in\IC$ by $\g_\zeta\subset\g$ and the subspace $[Z,\g]\subset\g$
of `off--diagonal' elements by $\god$. Thus, $\g^Z=\g_0$ and
\begin{equation}\label{eq:god}
\god=\sum_{\zeta\in\IC^*}\g_\zeta
\end{equation}

\subsection{} 

The assumptions above differ from those found in the literature in the
following ways:
\begin{enumerate}
\item with regard to $(Z)$, the assumption usually made is that $Z$ is
a {\it regular} semisimple element that is, if $G=GL_n(\IC)$, that $Z$
is diagonalisable with distinct eigenvalues. As pointed out in the Introduction,
this latter assumption is not stable under passage to a representation
and therefore ill--suited to the Tannakian methods employed in Section
\ref{se:Tannaka}.
\item On the other hand, the assumption $(f)$ is unduly restrictive.
A more natural assumption would be to consider the projection $f_0$
of $f$ onto $\g^Z$ and to require that the eigenvalues of $\ad(f_0)$
on $\g^Z$ are not positive integers. We impose the condition $(f)$
however since it simplifies the form of the Stokes data and holds in
the context of stability conditions considered in \cite{BTL}.
\end{enumerate}

\subsection{Stokes rays and sectors}

\begin{defn}\label{de:Stokes rays}
A {\it ray} is a subset of $\IC^*$ of the form $\R_{>0}\exp(i\pi\phi)$.
The {\it Stokes rays} of the connection $\nabla$ are the rays
$\IR_{>0}\zeta$, where $\zeta$ ranges over the non--zero eigenvalues
of $\ad(Z)$. The {\it Stokes sectors} are the open regions of $\C^*$
bounded by them. A ray is called {\it admissible} if it is not a Stokes ray.
\end{defn}

\begin{remark}
If $G$ is reductive, the eigenvalues of $\ad(Z)$ are invariant under
multiplication by $-1$. The Stokes rays therefore come in pairs and
the Stokes sectors are convex in this case. This need not be the case
for an arbitrary algebraic group.
\end{remark}

\subsection{Canonical fundamental solutions}

The Stokes data of the connection $\nabla$ are defined using fundamental
solutions with prescribed asymptotics. We first recall how these are
characterised.

Given a ray $r$ in $\C$, we denote by $\IH_r$ the corresponding half--plane
\begin{equation}\label{eq:halfplane}
\IH_r=\{z=uv:u\in r, \Re(v)>0\}\subset \C.
\end{equation}
The following basic result is well--known for $G=\GL_n(\IC)$ and $Z$
regular (see, \eg \cite[pp. 58--61]{Wasow}) and was extended in \cite
{Boalch2} to the case of complex reductive groups.\footnote{The
references \cite{Wasow} and \cite{Boalch2} cover however the more
general case when $\nabla$ is only defined on a disk around $t=0$
and has an arbitrary order pole at $t=0$.} It will be proved in Section
\ref{se:Tannaka}.

\begin{thm}
\label{jurk}
Given an admissible ray $r$, there is a unique holomorphic function
$Y_r:\IH_r\to G$ such that
\begin{gather}
\frac{dY_r}{dt}=\left(\frac{Z}{t^2}+\frac{f}{t}\right)Y_r
\label{eq:diff equ}\\
Y_r\cdot e^{Z/t}\to 1\quad\text{as}\quad\text{$t\to 0$ in $\IH_r$}
\label{eq:asym}
\end{gather}
\end{thm}

\begin{remark}
The function $Y_r\cdot e^{Z/t}$ possesses in fact an asymptotic expansion
in $\IH_r$ with constant term the identity, but we shall not need this stronger
property.
\end{remark}

\subsection{} 

The uniqueness statement of Theorem \ref{jurk} and the definition
of the Stokes data rely upon the following result which will be proved
in Section \ref{se:Phil lemma} (see \cite[Lemma 22]{Boalch2} for the
case of $G$ reductive and $Z$ regular).

\begin{prop}\label{pr:spectral}
Let $r,r'$ be two rays such that $r\neq -r'$, and $g\in G$ an element
such that
$$e^{-Z/t} \cdot g\cdot e^{Z/t}\to 1\text{ as }t\to 0\text{ in }\IH_r\cap \IH_{r'}.$$
Then, $g$ is unipotent and $X=\log(g)$ lies in
$$\bigoplus_{\zeta\in\ol{\Sigma}(r,r')}\g_\zeta\subset\g,$$
where $\ol{\Sigma}(r,r')\subset\IC^*$ is the closed convex sector bounded
by $r$ and $r'$.
\end{prop}

\noindent
Proposition \ref{pr:spectral} implies in particular that if the rays $r,r'$
are admissible and such that the sector $\ol{\Sigma}(r,r')$ does not
contain any Stokes rays of $\nabla$, the element $g\in G$ determined by
$$Y_r(t)=Y_{r'}(t)\cdot g\quad\text{for $t\in\IH_r\cap\IH_{r'}$}$$
is equal to 1. It follows in particular that, given a Stokes sector $\Sigma$
and a convex subsector $\Sigma'\subset\Sigma$, the solutions $Y_r$, as
$r$ varies in $\Sigma'$, patch to a fundamental solution $Y_{\Sigma'}$ of
\eqref{eq:diff equ} possessing the asymptotic property \eqref{eq:asym} in
the {\it supersector}
$$\wh{\Sigma'}=
\{uv:\thinspace u\in \Sigma', \Re(v)>0\}=
\bigcup_{r\subset\Sigma'}\IH_r.$$

\subsection{Stokes factors}\label{ss:Stokes factors}

Assume now that $\ell$ is a Stokes ray. Let $r_\pm$ be small clockwise
(resp. anticlockwise) perturbations of $\ell$ such that the convex sector
$\ol{\Sigma}(r_-,r_+)$ does not contain any Stokes rays of $\nabla$
other than $\ell$. 

\begin{defn}\label{de:Stokes factor}
The {\it Stokes factor} $S_\ell$ corresponding to $\ell$ is the element of
$G$ defined by
$$Y_{r_+}(t)= Y_{r_-}(t)\cdot S_\ell\text{ for }t\in \IH_{r_+}\cap \IH_{r_-}.$$
\end{defn}

\noindent By Proposition \ref{pr:spectral}, the definition of $S_\ell$ is
independent of the choice of $r_\pm$. Moreover, $S_\ell$ is unipotent
and $\log(S_\ell)\in\bigoplus_{\zeta\in\ell}\g_\zeta$.

\subsection{Stokes multipliers}\label{ss:multipliers}

An alternative but closely related system of invariants are the Stokes
multipliers of the connection $\nabla$. These depend upon the choice
of a ray $r$ such that both $r$ and $-r$ are admissible.

\begin{defn}
The {\it Stokes multipliers} of $\nabla$ corresponding to $r$ are the
elements $S_\pm\in G$ defined by
\[ Y_{r,\pm}(t)=Y_{-r}(t)\cdot S_\pm,\quad t\in\IH_{-r}\]
where $Y_{r,+}$ and $Y_{r,-}$ are the analytic continuations of $Y_r$
to $\IH_{-r}$ in the anticlockwise and clockwise directions respectively.
\end{defn}

\noindent By Proposition \ref{pr:spectral}, the multipliers $S_\pm$ remain
constant under a perturbation of $r$ so long as $r$ and $-r$ do not cross
any Stokes rays.

\subsection{}

To relate the Stokes factors and multipliers, set $r=\R_{>0}\exp(i\pi\theta)$
and label the Stokes rays as $\ell_j=\R_{>0}\exp(i\pi\phi_j)$, with $j=1,\ldots,
m_1+m_2$, where
$$\theta<\phi_1<\cdots<\phi_{m_1}< \theta+1
<\phi_{m_1+1}<\cdots<\phi_{m_1+m_2}<\theta+2$$
The following result is immediate upon drawing a picture
\begin{lemma}
\label{veryeasy}
The following holds
$$S_+=S_{\ell_{m_1}}  \cdots S_{\ell_1}\aand
S_-=S_{\ell_{m_1+1}}^{-1}\cdots S_{\ell_{m_1+m_2}}^{-1}$$
\end{lemma}

The Stokes factors therefore determine the Stokes multipliers for any ray
$r$. Conversely, the Stokes multipliers for a single ray $r$ determine all
the Stokes factors. This may be proved along the lines of \cite [Lemma 2]
{BJL}, \cite{Boalch2} by noticing that $S_\pm$ lie in the unipotent subgroups
$N_\pm\subset G$ with Lie algebras $\bigoplus_{\zeta\in\IH_{\pm ir}}\g_\zeta$
and that these groups may uniquely be written as products of the unipotent
subgroups $\exp(\bigoplus_{\zeta\in\ell}\g_\zeta)$ corresponding to the Stokes
rays $\ell$ in $\IH_{\pm ir}$. We shall instead give explicit formulae expressing
the Stokes factors in terms of $S_\pm$ in Proposition \ref{pr:Stokes HN Reineke}
below.

\subsection{Choice of a torus}

It will be convenient in the sequel to choose a torus $H\subset G$ whose
Lie algebra $\h$ contains $Z$. Let $\bfX(H)=\Hom_{\IZ}(H,\IC^*)$ be the
group of characters of $H$ and $\bfX(H)\cong\Lambda\subset\h^*$ the
lattice spanned by the differentials of elements in $\bfX(H)$. For any
$\lambda\in\Lambda$, we denote the unique element of $\bfX(H)$ with
differential $\lambda$ by $e^{\lambda}$. Decompose $\g$ as
\begin{equation}\label{eq:root decomposition 1}
\g=\g^\h\oplus[\h,\g]=
\g_0\oplus\bigoplus_{\alpha\in\Phi}\g_\alpha
\end{equation}
where $\Phi=\Phi(G,H)\subset\Lambda\setminus\{0\}$ is a finite set and
$H$ acts on $\g_\alpha$ via the character $e^\alpha$ so that, in particular
$\h\subset\g_0$. We refer to the elements of $\Phi$ as the {\it roots of $G$}
relative to $H$. We note that if $H$ is a maximal torus, the set of roots
$\Phi(G,H)$ is independent of the choice of $H$, but we shall not need
to assume this.

Set
\begin{equation}\label{eq:PhiZ}
\PhiZ=\{\alpha\in\Phi|\medspace Z(\alpha)\neq 0\}
\end{equation}
so that
\begin{equation}\label{eq:Z and h}
[Z,\g]=\bigoplus_{\alpha\in\PhiZ}\g_\alpha
\aand
\g^Z=\g_0\oplus\bigoplus_{\alpha\in\Phi\setminus\PhiZ}\g_\alpha
\end{equation}

\begin{remark} 
If $f=\sum_{\zeta\in\IC}f_\zeta$ is the eigenvector decomposition
of $f\in\g$ \wrt $\ad(Z)$, then
\begin{equation}\label{eq:eliminate H}
f_\zeta=\sum_{\alpha:Z(\alpha)=\zeta}f_\alpha
\end{equation}
\end{remark}

\section{Isomonodromic deformations}\label{se:imd}

We discus in this section isomonodromic deformations of the
connection $\nabla$, that is families of connections of the form
\eqref{eq:reference} where $Z$ and $f$ vary in such a way that
the Stokes data remain constant. These deformations will be
used in Section \ref{ss:Jn} to establish the analytic properties 
of the local inverse of the Stokes map.

\subsection{Variations of $Z$} 

We wish to vary $Z$ among semisimple elements of $\g$ in such
a way that the decomposition \eqref{eq:Z decomposition} remains
constant.

This is readily seen to be the case if $Z$ varies among the regular
elements of the Lie algebra $\h$ of a torus $H\subset G$. Indeed,
let $\Phi\subset\h^*$ be the set of roots of $H$ and set
$$\hreg=\h\setminus\bigcup_{\alpha\in\Phi}\Ker(\alpha)$$
If $Z\in\hreg$, $\Phi^Z$ is equal to $\Phi(G;H)$ so that the
decomposition \eqref{eq:Z decomposition} remains constant
as $Z$ varies in $\hreg$ by \eqref{eq:Z and h}.

Conversely, the following holds

\begin{prop}
Let $Z\in\g$ be a fixed semisimple element.
\begin{enumerate}
\item The set of semisimple elements $Z'\in\g$ which give rise
to the same decomposition \eqref{eq:Z decomposition} as $Z$
is the set of regular elements in the Lie algebra $\h$ of a torus
$H\subset G$ such that $Z\in\h$.
\item $H$ is the unique torus in $G$ which is maximal for the
property that $Z\in\hreg$.
\end{enumerate}
\end{prop}
\begin{pf}
(i)
Let $\g^Z\subseteq\g$ be the centraliser of $Z$, $\Z(\g^Z)$ its
centre and $\h\subseteq\Z(\g^Z)$ the subspace of its semisimple
elements. $\h$ is the Lie algebra of a torus $H\subset G$ obtained
as follows. Let $G^Z\subseteq G$ be the centraliser of $Z$, $G
^Z_1\subseteq G$ its identity component and $\Z(G^Z_1)_1$ the
identity component of the centre of $G^Z_1$. This is a connected,
commutative algebraic group. Its subgroup $H$ of semisimple
elements is therefore a torus with Lie algebra $\h$.

Note next that $\g^Z=\g^\h$ and $[Z,\g]=[\h,\g]$ so that $Z\in\hreg$.
Indeed, since $Z\in\h$, we have $\g^\h\subseteq\g^Z$ and $[\h,\g]
\supseteq[Z,\g]$. By definition however, $\h\subset\Z(\g^Z)$ so that
$\g^Z\subseteq\g^\h$ and both of the previous inclusions are equalities. 

Let now $Z'\in\g$ be such that $\g^{Z'}=\g^Z$, then $Z'\in\Z(\g
^{Z'})=\Z(\g^Z)$. Since $Z'$ is semisimple, $Z'\in\h$. If suffices
now to notice that, for any $Z'\in\h$ one has $\g^{Z'}\supseteq
\g^\h=\g^Z$ and $[Z',\g]\subseteq [\h,\g]=[Z,\g]$ with equalities
if, and only if, $\alpha(Z')\neq 0$ for all $\alpha\in\Phi$.

(ii) If $\wt{H}\subset G$ is a torus such that $Z\in\wt{\h}_\reg$,
(i) implies that $\wt{\h}_\reg\subset\hreg$ and therefore that $
\wt{H}\subset H$.
\end{pf}

\subsection{Isomonodromic families of connections}\label{ss:IM}

Fix henceforth a torus $H\subset G$. Let $P$ be the holomorphically
trivial principal $G$--bundle over $\IP^1$ and let $\U \subset\hreg$
be an open set. Consider a family of connections on $P$ of the form
\eqref{nab}, namely
$$\nabla(Z)=d-\bigg(\frac{Z}{t^2}+\frac{f(Z)}{t}\bigg)dt$$
where $Z$ varies in $\U$ and the dependence of $f(Z)\in\god$
\wrt $Z$ is arbitrary.

\begin{defn}\label{de:IMD multipliers}
The family of connections $\nabla(Z)$ is {\it isomonodromic} if for any
$Z_0\in\U$, there exists a neighborhood $Z_0\in\U_0\subset\U$ and
a ray $r$ such that $\pm r$ are admissible for all $\nabla(Z)$, $Z\in\U
_0$ and the Stokes multipliers $S_\pm(Z)$ of $\nabla(Z)$ relative to
$r$ are constant on $\U_0$.
\end{defn}

The isomonodromy of the family $\nabla(Z)$ may also be defined as the
constancy of the Stokes factors. This requires a little more care since,
as pointed out in \cite[pg. 190]{Boalch} for example, Stokes rays may
split into distinct rays under arbitrarily small deformations of $Z$. Call
a sector $\Sigma\subset\IC^*$ {\it admissible} if its boundary rays are
admissible. 

\begin{prop}\label{pr:IMD factors}
The family of connections $\nabla(Z)$ is isomonodromic if, and only
if, for any connected open subset $\U_0\subset\U$ and any convex
sector $\Sigma$ which is admissible for all $\nabla(Z)$, $Z\in\U_0$,
the clockwise product
$$\clockwise_{\ell\subset\Sigma}S_\ell(Z)$$
of Stokes factors corresponding to the Stokes rays contained in
$\Sigma$ is constant on $\U_0$.
\end{prop}
\begin{pf}
This follows from the fact that Stokes factors and multipliers determine
each other by Lemma \ref{veryeasy} and Proposition \ref{pr:Stokes HN Reineke}.
\end{pf}

\subsection{Isomonodromy equations}

The following characterisation of isomonodromic deformations was
obtained by Jimbo--Miwa--Ueno \cite{JMU} when $G=GL_n(\IC)$
and $H$ is a maximal torus and adapted to the case of a complex,
reductive group by Boalch \cite[Appendix]{Boalch}. The proof carries
over verbatim to the case of an arbitrary algebraic group $G$ and
torus $H\subset G$.

\begin{thm}\label{th:IMD}
Assume that $f$ varies holomorphically in $Z$. Then, the family of
connections $\nabla(Z)$ is isomonodromic if, and only if $f$ satisfies
the PDE
\begin{equation}
\label{iso}
d f_\alpha=\sum_{\substack{\beta,\gamma\in\Phi:\\\beta+\gamma=\alpha}}
[f_\beta, f_\gamma]\, d\log\gamma.
\end{equation}
\end{thm}

\remark\label{rk:IMD} The equations \eqref{iso} form a first order system
of {\it integrable} non--linear PDEs and therefore possess a unique
holomorphic solution $f(Z)$ defined in a neighboorhood of a fixed 
$Z_0\in\hreg$ and subject to the initial condition $f(Z_0)=f_0\in[Z_0,\g]$.

\remark Jimbo--Miwa--Ueno and Boalch also give an alternative characterisation
of isomonodromy in terms of the existence of a flat connection on $\IP^1\times\U$
which has a logarithmic singularity on the divisor $\{t=\infty\}$ and a pole of order
2 on $\{t=0\}$, and restricts to $\nabla(Z)$ on each fibre $\{Z\}\times\IP^1$. This
connection is given by
\[\ol{\nabla} = d - \bigg[\bigg(\frac{Z}{t^2}+\frac{f}{t}\bigg)dt
 +  \sum_{\alpha\in \Phi} f_\alpha\frac{d\alpha}{\alpha}
+\frac{ dZ}{t}\bigg].\]
One can check directly that the flatness of this connection
is equivalent to \eqref{iso}.

\section{The Stokes map}\label{se:Stokes map}

In this section, we express the logarithms of the Stokes factors of
the connection $\nabla$ as explicit, universal Lie series in the
variables $f_\alpha$. Using the results of Section \ref{comb}, we
then show how to invert these series to express the $f_\alpha$ as
Lie series in the logarithms of the Stokes factors, thus explicitly
solving a Riemann--Hilbert problem.

\subsection{Completion with respect to \fd representations}

Our formulae for the Stokes factors $\nabla$ are more conveniently
expressed inside the completion $\wh{U\g}$ of $U\g$ \wrt the \fd
representations of $G$. We review below the definition of $\wh{U\g}$.

Let $\Vec$ be the category of \fd complex vector spaces and $\RRep
(G)$ that of \fd representations of $G$. Consider the forgetful functor
$$F:\RRep(G)\rightarrow\Vec.$$
By definition, $\wh{U\g}$ is the algebra of endomorphisms of $F$.
Concretely, an element of $\wh{U\g}$ is a collection $\Theta=\{\Theta
_V\}$, with $\Theta_V\in\End_{\IC}(V)$ for any $V\in\RRep(G)$, such
that for any $U,V\in\RRep(G)$ and $T\in\Hom_G(U,V)$, the following
holds
$$\Theta_V\circ T=T\circ\Theta_U$$

There are natural homomorphisms $G\to\wh{U\g}$ and $U\g\to\wh{U\g}$ 
mapping $g\in G$ and $x\in U\g$ to the elements $\Theta(g)$, $\Theta
(x)$ which act on a \fd representation $\rho:G\to GL(V)$ as $\rho(g)$
and $\rho(x)$ respectively. The following is well known.

\begin{lemma}\label{le:injective}
The homomorphisms $G\to\wh{U\g}$ and $U\g\to\wh{U\g}$ are injective.
\end{lemma}
\begin{pf}
The first claim follows immediately from the fact that $G$ has a faithful
\fd representation. For the second, we use the fact that $U\g$ acts faithfully
on $\IC[G]$ by left--invariant differential operators. Since this action
decreases the degree of polynomials, any $f\in\IC[G]$ is contained in a
\fd $G$--module and the claim follows.
\end{pf}

We will use the homomorphisms above to think of $U\g$ as a subalgebra
of $\wh{U\g}$ and $G$ as a subgroup of the group of invertible elements
of $\wh{U\g}$ respectively.

\subsection{Representing Stokes factors}\label{rep}

Fix a Stokes ray $\ell$. We show below how to represent the
corresponding Stokes factor $S_\ell$ in two different ways: by
elements $\epsilon_\alpha\in\god$ and by elements $\delta_
\gamma\in U\g$.

Consider the subalgebra
$$\n_\ell=\bigoplus_{\alpha:Z(\alpha)\in\ell} \g_\alpha\subset \g.$$
The elements of $\n_\ell$ are nilpotent, that is they act by nilpotent
endomorphisms on any \fd representation of $G$. It follows that the
exponential map $\exp\colon\n_\ell\to G$ is an isomorphism onto the
unipotent subgroup $N_\ell=\exp(\n_\ell)\subset G$.

By Proposition \ref{pr:spectral}, the Stokes factor $S_\ell$ lies in
$N_\ell$. For the first representation of $S_\ell$, write
\begin{equation}
\label{bubbly}
S_\ell=\exp\bigg(\sum_{\alpha:Z(\alpha)\in\ell}\epsilon_\alpha\bigg)
\end{equation}
for uniquely defined elements $\epsilon_\alpha\in\g_\alpha$.
For the second, we compute the exponential
\eqref{bubbly} in $\wh{U\g}$ and decompose the result along the weight
spaces
$$\wh{U\g}_\gamma=
\{x\in\wh{U\g}|\ad(h)x=\gamma(h)x,\;\forall h\in\h\},\;\gamma\in\h^*$$
of the adjoint action of $\h$. This yields elements $\delta_\gamma\in
(U\n_\ell)_\gamma$ such that
\begin{equation}
\label{hold}
S_{\ell}
=1+\sum_{\gamma\in\LambdaZ:Z(\gamma)\in\ell}\delta_\gamma,
\end{equation}
where $\LambdaZ\subset\h^*$ is the lattice generated by $\PhiZ$ and
the above identity is to be understood as holding in any \fd \rep of $G$
(where the \rhs is necessarily finite).

These two representations of $S_\ell$ are related as follows.

\begin{lemma}\label{le:epsilon delta}\hfill
\begin{enumerate}
\item
Let $\gamma\in\LambdaZ$ be such that $Z(\gamma)$ lies on the
Stokes ray $\ell$. Then, $\delta_\gamma$ is given by the finite sum
\begin{equation}
\label{exp}
\delta_\gamma=
\sum_{n\geq 1}
\sum_{\substack{\alpha_i\in\PhiZ,
\\Z(\alpha_i)\in\ell,\\
\alpha_1+\cdots+\alpha_n=\gamma}}
\frac{1}{n!}\;\epsilon_{\alpha_1} \cdots \epsilon_{\alpha_n}.
\end{equation}
\item
Conversely, let $\alpha\in\PhiZ$ be such that $Z(\alpha)\in\ell$.
Then, $\epsilon_\alpha$ is given by the finite sum
\begin{equation}
\label{log}
\epsilon_\alpha=
\sum_{n\geq 1}
\sum_{\substack{\gamma_i\in\LambdaZ,\\
Z(\gamma_i)\in\ell,\\
\gamma_1+\cdots+\gamma_n=\alpha}}
\frac{(-1)}{n}^{n-1}\,\delta_{\gamma_1}\cdots\delta_{\gamma_n}.
\end{equation}
\end{enumerate}
\end{lemma}
\begin{pf}
These are the standard expansions of $\exp\colon\n_\ell\to N_\ell$ and $\log\colon
N_\ell\to\n_\ell$.
\end{pf}

\subsection{The Stokes map}\label{ss:Stokes}

Since the subsets $\{\alpha\in\Phi:Z(\alpha)\in\ell\}$ partition $\Phi$ as
$\ell$ ranges over the Stokes rays of $\nabla$, we may assemble the
elements $\epsilon_\alpha$ corresponding to different Stokes rays and
form the sum
\begin{equation}\label{eq:g od}
\epsilon=
\sum_{\alpha\in \PhiZ}\epsilon_\alpha\in\bigoplus_{\alpha\in\PhiZ}\g_\alpha.
\end{equation}
We shall refer to the map
\begin{equation}\label{eq:Stokes map}
\calS:\bigoplus_{\alpha\in\PhiZ}\g_\alpha\longrightarrow
\bigoplus_{\alpha\in\PhiZ}\g_\alpha
\end{equation}
mapping $f$ to $\epsilon$ as the {\it Stokes map}. Note that $\calS$ depends
upon $Z$.

\subsection{The functions $M_n$}

We give below an explicit formula for the Stokes factors of the connection
$\nabla$ in terms of iterated integrals. The definition and elementary
properties of iterated integrals are reviewed in Section \ref{polylog}.

\begin{defn}\label{de:M}
Set $M_1(z_1)=2\pi i$ and, for $n\geq 2$, define the function $M_n:
(\IC^*)^n\to\IC$ by the iterated integral
$$M_n(z_1,\ldots,z_n)=
2\pi i \big.\int_{C}\frac{dt}{t-s_1}\circ \cdots \circ \frac{dt}{t-s_{n-1}},$$
where $s_i=z_1+ \cdots +z_i$, $1\leq i\leq n$ and the path of integration
$C$ is the line segment $(0,s_n)$, perturbed if necessary to avoid any
point $s_i\in[0,s_n]$ by small clockwise arcs.
\end{defn}

\begin{remark}
The iterated integrals defining the functions $M_n$ are convergent
by Lemma \ref{wellbehaved} below since the assumption that $z_i
\in\IC^*$ for all $i$ implies in particular that $s_1\neq 0$ and $s_{n-1}
\neq s_n$.
\end{remark}

\subsection{Formula for the Stokes factors}\label{ss:formula for factors}

The following result will be proved in Section \ref{se:Tannaka}.

\begin{thm}\label{one2}
The Stokes factor $S_\ell$ corresponding to the ray $\ell$ is given by
\begin{equation}\label{eq:S ell}
S_\ell=
1+\sum_{n\geq 1}\,
\sum_{\substack{\alpha_1,\ldots,\alpha_n\in\PhiZ\\
Z(\alpha_1+\cdots+\alpha_n)\in\ell}}
M_n(Z(\alpha_1),\ldots, Z(\alpha_n))\,
f_{\alpha_1}\cdots f_{\alpha_n}
\end{equation}
where the equality is understood as holding in any \fd \rep of $G$ and
the sum over $n$ is absolutely convergent.
\end{thm}

\noindent In terms of weight components, \eqref{eq:S ell} reads
\begin{equation}
\label{mainy}
\delta_\gamma=
\sum_{n\geq 1}
\sum_{\substack{\alpha_1,\ldots,\alpha_n\in\PhiZ\\
\alpha_1+\cdots+\alpha_n=\gamma}}
M_n(Z(\alpha_1),\ldots, Z(\alpha_n))\,
f_{\alpha_1} f_{\alpha_2}\cdots f_{\alpha_n}.
\end{equation}
for any $\gamma\in\LambdaZ$ such that $Z(\gamma)\in\ell$.

\remark Theorem \ref{one2} shows that the Stokes factors of $\nabla$ are given
by periods. Their appearence in this context stems from the fact that their computation
reduces, via the results of Sections \ref{se:Fourier-Laplace}--\ref{se:Stokes factors},
to one of partial monodromies of the \FL transform $\nablah$ of $\nabla$, which
are well--known to be given by iterated integrals. It seems an interesting problem
to determine whether the Stokes factors of a connection with arbitrary order poles
on $\IP^1$ are also expressible in terms of explicit periods.

\subsection{The functions $L_n$}

We next state a formula for the Stokes map giving the element $\epsilon$ in
terms of $f$. We first define  the special functions appearing in this formula.

\begin{defn}\label{de:Ln}
The function $L_n:(\IC^*)^n\to\IC$ is given by $L_1(z_1)=2\pi i$ and, for
$n\geq 2$,
$$L_n(z_1,\ldots,z_n)=
\sum_{k=1}^n\;\sum_{\substack{
0=i_0<\cdots<i_k= n\\
s_{i_j}-s_{i_{j-1}}\in\IR_{>0}\cdot s_n}}
\frac{(-1)}{k}^{k-1}\;
\prod_{j=0}^{k-1} M_{i_{j+1}-i_j}(z_{i_j +1} ,\ldots,z_{i_{j+1}}),$$
where $s_j=z_1+\cdots+z_j$.
\end{defn}

\remark Note that $L_1\equiv M_1$ and that on the open subset
\[(z_1,\ldots,z_n) \in (\C^*)^n\text{ such that } s_i\notin [0,s_n]\text
{ for }0<i<n\] 
the inner sum above is empty unless $k=1$ so that
$L_n(z_1,\ldots,z_n)=M_n(z_1,\ldots,z_n).$
Thus $L_n$ agrees with $M_n$ on the open subset where it is holomorphic
and differs from it by how it has been extended onto the cutlines.

The functions $L_n$ are more complicated to define than the functions
$M_n$. Unlike the latter however, they give rise to Lie series by Theorem
\ref{one} (i).

\remark For $n\geq 2$, the function $L_n$ satisfies
\begin{equation}\label{eq:L zero sum}
L_n(z_1,\ldots,z_n)=0
\quad\text{if}\quad
z_1+\cdots+z_n=0
\end{equation}
Indeed, the summation condition above becomes $z_{i_j}+\cdots+z_{i_{j-1}+1}
=0$ for any $j=1,\ldots,n$ and $M_m(w_1,\ldots,w_m)=0$ whenever $w_1+
\cdots+w_m=0$.

\subsection{Formula for the Stokes map}\label{ss:Stokes map}

\begin{thm}\label{one}\hfill
\begin{enumerate}
\item Let $x_1,\ldots, x_m$ be elements in a Lie algebra $\L$. For any $(z_1,
\ldots,z_m)\in (\C^*)^m$, the finite sum
\begin{equation}\label{eq:Lie poly}
\sum_{\sigma\in \Sym_m}
L_m(z_{\sigma(1)}, \cdots, z_{\sigma(m)})
x_{\sigma(1)} \cdots x_{\sigma(m)}
\end{equation}
is a Lie polynomial in $x_1,\ldots,x_m$ and therefore lies in $\L\subset U\L$.
\item The element $\epsilon=\calS(f)$ is given by the following Lie series in
the variables $\{f_\alpha\}_{\alpha\in\PhiZ}$
\begin{equation}
\label{main}
\epsilon_\alpha=
\sum_{n\geq 1}\sum_{\stackrel{\alpha_i\in\PhiZ}{\alpha_1+\cdots+\alpha_n=\alpha}}
L_n(Z(\alpha_1),\ldots, Z(\alpha_n)) f_{\alpha_1} f_{\alpha_2}\cdots f_{\alpha_n}
\end{equation}
where the sum over $n$ is absolutely convergent.
\end{enumerate}
\end{thm}
\begin{pf}
Substituting \eqref{mainy} into \eqref{log} shows that \eqref{main} holds
in $\wh{U\g}$, where the sum over $n$ is absolutely convergent in any
\fd representation of $G$.

To prove (i), we may assume that $z_1+\cdots+z_m\neq 0$ since $L_m
(z_1,\ldots,z_m)=0$ if $z_1+\cdots+z_m=0$. We need to show that if $\L$
is the free Lie algebra generated by $x_1,\ldots,x_m$, the sum \eqref
{eq:Lie poly} lies in $\L$. $\L$ is $\IZ_{\geq 0}$--graded by setting $\deg
(x_i)=1$ and it suffices to prove this for the \fd quotient $\L_{\leq m}$ of
$\L$ by elements of degree $\geq m+1$.

We shall do so by applying the version of (ii) just obtained to the semi--direct
product $B=T\ltimes L_{\leq m}$ of the algebraic torus $T=(\IC^*)^m$
by the unipotent group $L_{\leq m}$ with Lie algebra $\L_{\leq m}$, where
$T$ acts on $\L_{\leq m}$ by $(w_1,\ldots,w_m)\,x_i=w_ix_i$.
Let $\t=\IC^m$ be the Lie algebra of $T$, $\{E_i\}_{i=1}^m$ its
canonical basis and $\{\theta_i\}$ the corresponding dual basis of $\t^*$.
The root system $\Phi(B,T)\subset\t^*$ is contained in the set of elements
of the form $\sum_i p_i\theta_i$, with $p_i\in\IZ_{\geq 0}$ and $\sum_i p_i
\leq m$.

Consider the connection $\nabla$ with structure group $B$ and coefficients
$Z\in\t$ and $f\in\L_{\leq m}$ given by
$$Z=z_1 E_1+\cdots+z_m E_m\aand
f=x_1+\cdots+x_m.$$
The logarithms of the Stokes factors of $\nabla$ determine an element
$\epsilon\in\L_{\leq m}$ whose components are given by \eqref{main},
understood as an identity in $\wh{U\L_{\leq m}}$. Since in this case the
sum over $n$ is finite, \eqref{main} actually holds as an identity in $U\L
_{\leq m}$.

Note that the roots $\theta_1,\ldots,\theta_m$ and $\theta_1+\cdots+\theta
_m$ lie in $\PhiZ$ since, by assumption $Z(\theta_i)=z_i\neq 0$ and
$Z(\theta_1+\cdots+\theta_m)=z_1+\cdots+z_m\neq 0$. Since the only
non--zero components of $f$ are $f_{\theta_i}=x_i$, $i=1,\ldots,m$,
the component $\epsilon_{\theta_1+\cdots+\theta_m}$ is precisely
equal to \eqref{eq:Lie poly}. This proves (i), and thus that \eqref{main}
holds as an identity in $\g\subset\wh{U\g}$.
\end{pf}

\remark\label{rk:global main}
The \rhs of \eqref{mainy} makes sense for any $\gamma\in\LambdaZ$.
It is easy to show using Theorem \ref{one} that each summand in $n$ is equal
to zero unless $\gamma$ is such that $Z(\gamma)$ lies on a Stokes ray of
$\nabla$. Thus, the identity \eqref{mainy} holds for any $\gamma\in\LambdaZ$.

Similarly, \eqref{main} holds for any $\alpha\in\Phi$ since, for $\alpha\in\Phi
\sqcup\{0\}\setminus\PhiZ$, the \lhs is equal to zero by definition and the \rhs
vanishes by \eqref{eq:L zero sum}. Thus, \eqref{main} may equivalently
be written as
$$\epsilon=
\sum_{n\geq 1}\sum_{\alpha_1,\ldots,\alpha_n\in\Phi}
L_n(Z(\alpha_1),\ldots, Z(\alpha_n)) f_{\alpha_1} f_{\alpha_2}\cdots f_{\alpha_n}
$$ 


\subsection{Inverse of  the Stokes map}

By Theorem \ref{one} and \cite[Thm. 1.5.6]{Rudin}, the series \eqref{main}
converges uniformly on compact subsects of $\god=\bigoplus_{\alpha\in\PhiZ}
\g_\alpha$. Thus, the Stokes map $\calS:\god\to\god$ is holomorphic, satisfies
$\calS(0)=0$ and its differential at $f=0$ is invertible since $L_1$ is identically
equal to $2\pi i$. By the inverse function theorem, $\calS$ possesses an analytic
inverse $\calS^{-1}$ defined on a neighborhood of $\epsilon=0$. 

\begin{thm}\label{th:Stokes inverse}
The Taylor series of $\calS^{-1}$ at $\epsilon=0$ is given by a Lie series
in the variables $\{\epsilon_\alpha\}_{\alpha\in\PhiZ}$ of the form
\begin{equation}
\label{main2}
f_\alpha=
\sum_{n\geq 1}
\sum_{\stackrel{\alpha_i\in\PhiZ}{\alpha_1+\cdots+\alpha_n=\alpha}}
J_n(Z(\alpha_1),\ldots, Z(\alpha_n))\,
\epsilon_{\alpha_1} \epsilon_{\alpha_2}\cdots \epsilon_{\alpha_n}
\end{equation}
where the functions $J_n\colon (\C^*)^n \to \C$ are independent of $\g$
and such that $J_1\equiv 1/2\pi i$ and, for $n\geq 2$,
\begin{equation}\label{eq:vanishing}
J_n(z_1,\ldots,z_n)=0\quad\text{if}\quad z_1+\cdots+z_n=0
\end{equation}
\end{thm}

\begin{remark} Analogously to Remark \ref{rk:global main}, \eqref{main2}
may equivalently be written as
$$f=
\sum_{n\geq 1}
\sum_{\alpha_1,\ldots,\alpha_n\in\Phi}
J_n(Z(\alpha_1),\ldots, Z(\alpha_n))\,
\epsilon_{\alpha_1} \epsilon_{\alpha_2}\cdots \epsilon_{\alpha_n}$$
\end{remark}

\subsection{The functions $\mathbf{J_n}$}\label{ss:Jn}

Theorem \ref{th:Stokes inverse} will be proved in Section \ref{comb}
by formally inverting the power series \eqref{main}. This yields an
explicit definition of the functions $J_n$ as sums of products of the
functions $L_n$ indexed by plane rooted trees. For example,
\begin{multline*}
(2\pi i)^3 J_3(z_1,z_2,z_3)= \\
L_2(z_1+z_2,z_3)L_2(z_1,z_2)-
L_3(z_1,z_2,z_3)+
L_2(z_1,z_2+z_3) L_2(z_2,z_3)
\end{multline*}
corresponding to the three distinct plane rooted trees with 3 leaves.

The following establishes the analytic properties of the functions $J
_n$. These are not readily apparent from the combinatorial definition
of these functions, and will be obtained instead by using the
isomonodromic deformations considered in Section \ref{se:imd}.

\begin{thm}
\label{js}
The function $J_n\colon (\C^*)^n\to \C$ is continuous and holomorphic
on the complement of the hyperplanes
\begin{equation*}
H_{ij}=\{z_i+\cdots+z_j=0\},\quad 1\leq i<j\leq n
\end{equation*}
in the domain
\begin{equation}\label{eq:Joyce domain}
\D_n=\{(z_1,\ldots,z_n)\in (\C^*)^n|\medspace 
z_i/z_{i+1} \notin \R_{>0}
\text{ for } 1\leq i< n\}
\end{equation}
Moreover, it satisfies the differential equation
\begin{equation}
\label{joycediff}
d J_n(z_1,\ldots,z_n)=
\sum_{i=1}^{n-1} J_{i}(z_1,\ldots, z_{i}) J_{n-i}(z_{i+1},\ldots, z_n)
d\log \bigg(\frac{z_{i+1}+\cdots+z_n}{z_1+\cdots +z_i} \bigg)
\end{equation}
\end{thm}
\begin{pf}
The stated properties of the functions $J_n$ will be obtained by
applying Theorem \ref{th:Stokes inverse} to the group $B\subset
GL_{n+1}(\IC)$ of upper triangular matrices.
Choose as torus $H\subset B$ the subgroup of diagonal matrices,
let $N\subset B$ be the subgroup of strictly upper triangular
matrices and denote their Lie algebras by $\b,\h$ and $\n$
respectively. Let $\{e_i\}_{i=1}^{n+1}$ be the canonical basis
of $\IC^{n+1}$, $E_{ij}\,e_k=\delta_{jk}e_i$ the corresponding
elementary matrices, and $\{\theta_i\}_{i=1}^{n+1}$ the basis
of $\h^*$ given by $\theta_i(E_{jj})=\delta_{ij}$. The root system
$\Phi(B,H)$ consists of the linear forms $\alpha=\theta_i-\theta_j$,
$1\leq i<j\leq n+1$. Note that in this case the sum \eqref{main2}
is finite and therefore defines, for any fixed $Z\in\h$, a global
inverse to the Stokes map $\calS:[Z,\g]\to[Z,\g]$.

Let $(z_1,\ldots,z_n)\in (\C^*)^n$, and set
$$Z=\diag(z_1+\cdots+z_n,z_2+\cdots+z_n,\ldots,z_n,0)\in\h$$
Note that since $Z(\theta_i-\theta_j)=z_i+\cdots+ z_{j-1}$, $Z$
lies in $\hreg$ if, and only if $(z_1,\ldots,z_n)\notin\bigcup_{i\neq j}
H_{ij}$. Set $\epsilon=\sum_{i=1}^n E_{i\,i+1}\in\n$. Since $z_i
\neq 0$, $\epsilon$ lies in
$$[Z,\b]=
\bigoplus_{\substack{1\leq i<j\leq n+1\\z_i+\cdots+z_{j-1}\neq 0}}
\IC E_{ij}$$
and therefore defines Stokes data for a unique $B$--connection
of the form \eqref{nab}. The corresponding Stokes rays are $\IR
_{>0}(z_i+\cdots+ z_{j-1})$ and include in particular the rays $\ell
_i=\IR_{>0}z_i$. Since $\epsilon_\alpha=0$ unless $\alpha=\theta
_i-\theta_{i+1}$, the only non--trivial Stokes factors correspond to
the rays $\ell_i$ and are given by
$$S_{\ell_i}=\exp\left(\sum_{j:z_j\in\ell_i}E_{j\,j+1}\right)$$
Assume now that $(z_1,\ldots,z_n)\in\D_n$, and set $S_i=\exp(E_
{i\,i+1})=1+E_{i\,i+1}$. Given that $[E_{i\,i+1},E_{j\,j+1}]=0$ if $|i-j|
\geq 2$, the above Stokes factors are given by $S_{\ell_i}=\prod_
{j:z_j\in\ell_i}S_j$ since $\ell_i\neq\ell_{i\pm 1}$ on $\D_n$. Thus,
varying $(z_1,\ldots,z_n)$ in $\D_n'=\D_n\setminus\bigcup_{i<j}
H_{ij}$, and keeping $\epsilon$ fixed yields an isomonodromic
family of connections.

Since the Stokes map $\calS:[Z,\g]\to[Z,\g]$ has a global inverse
for $G=B$, the corresponding $f=f(Z)$ varies holomorphically in
$Z\in\hreg$ and satisfies the isomonodromy equations \eqref{iso}.
Indeed, fix $Z_0\in\hreg$, let $Z_0\in\U\subset\hreg$ be a small
open neighborhood and $f':\U\to\n$ the local holomorphic
solution of the isomonodromy equations \eqref{iso} such that $f'(
Z_0)=f(Z_0)$ (see Remark \ref{rk:IMD}). By Theorem \ref{th:IMD},
the connection $d-(Z/t^2+f'(Z)/t)dt$ has constant Stokes data as
$Z$ varies in $\U$ so that $f'(Z)=f(Z)$ on $\U$ since $\calS$ is
injective.

Computing now \eqref{main2} in the vector representation yields
$$f_{\theta_i-\theta_j}=J_{j-i}(z_i,\ldots,z_{j-1})E_{i\,j}$$
for any $i<j$. The claimed regularity of the functions $J_i$,
$i\leq n$ on $\D_n'$ and the differential equation \eqref
{joycediff} now follow from that of $f$ and the equations
\eqref{iso}.
\end{pf}

\begin{remark}
By arguing as in \cite[Prop. 3.10]{Joyce}, it is easy to prove
by induction on $n$, using the PDE \eqref {joycediff} that the
function $J_n$ possesses a holomorphic extension $\wt{J}_
n$ to $\D_n$. Although we have not checked this in full, we
believe that this extension coincides with our combinatorial
definition of $J_n$.  
\end{remark}

\subsection{Irregular Riemann--Hilbert correspondence}

It is important to distinguish the Stokes map \eqref{eq:Stokes
map} from a related map also studied by \BJL \cite{BJL} and
Boalch \cite{Boalch3,Boalch2}. Rather than studying connections
of the form \eqref{nab}, one can consider instead meromorphic
connections on the trivial principal $G$--bundle over the unit
disc $D\subset \C$ which have the form
\begin{equation}
\label{nabb}
d-\bigg(\frac{Z}{t^2}+\frac{f(t)}{t}\bigg)dt
\end{equation}
where  $f\colon D \to \g$ is holomorphic.

One can define Stokes data for such a connection as in Section \ref
{ss:Stokes factors}, and consider the map sending the set of gauge
equivalence classes of such connections to that of possible Stokes
data. Boalch refers to this map as the {\it irregular Riemann--Hilbert
map}. Extending results of \cite{BJL} for $G=GL_n(\IC)$, he shows
that, for $G$ reductive and $Z$ regular semisimple, this map is in
fact an isomorphism.

In contrast, as shown in \cite{JLP}, the Stokes map $\calS$ is neither
injective nor surjective in general, even when $G=GL_2(\IC)$. Put another
way, not every connection of the form \eqref{nabb} can be put into the
constant--coefficient form \eqref{nab} by a gauge transformation, and
even when that is possible, the resulting connection \eqref{nab} is not
in general unique.

It follows from Theorem \ref{th:Stokes inverse} that whenever the sum
\eqref{main2} is absolutely convergent over $n$, it successfully inverts
the Stokes map $\calS$, in that the connection \eqref{nab} determined
by $f\in\god$ has Stokes factors given by \eqref{bubbly}. In spite of our
assumption that $f\in\god$, which does not hold in the counterexamples
of \cite{JLP}, we do not expect the Stokes map to be bijective and therefore
the sum to be absolutely convergent over $n$ for arbitrary $\epsilon$.


\section{Stokes multipliers}\label{se:multipliers}

Throughout this section, we fix a ray $r=\R_{>0}e^{i\pi\theta}$ such that
$\pm r$ are admissible for the connection \eqref{nab} and consider the
Stokes multipliers $S_\pm$ relative to $r$. We shall give an explicit
formula for $S_\pm$ analogous to that for the Stokes factors given
by Theorem \ref{one2}.

\subsection{Representing Stokes multipliers}\label{ss:rep kappa}

We first show how to represent $S_\pm$ by an element $\kappa\in\wh
{U\g}$. Let $\pm i\IH_r$ be the connected components of $\IC\setminus
\IR\,e^{i\pi\theta}$. These determine a partition of $\PhiZ=\PhiZ_+\sqcup
\PhiZ_-$ given by
$$\PhiZ_\pm=\{\alpha\in\PhiZ:Z(\alpha)\in\pm i\IH_r\}$$
Let $\LambdaZ_\pm\subset\h^*\setminus\{0\}$ be the cones spanned
by the linear combinations of elements in $\Phi_\pm^Z$ with coefficients
in $\IN_{>0}$. Similarly to \S \ref{rep}, it follows from Proposition \ref
{pr:spectral} that there is a unique element
$$\kappa=
\sum_{\gamma\in\LambdaZ_+\sqcup\LambdaZ_-}
\kappa_\gamma\in \wh{U\g}$$
such that the Stokes multipliers $S_\pm$ are respectively equal to
$$S_+=
1+\sum_{\gamma\in\LambdaZ_+}\kappa_\gamma
\aand
(S_-)^{-1}=1+\sum_{\gamma\in\LambdaZ_-}\kappa_\gamma,$$

\subsection{}

Given $\gamma\in\LambdaZ_+$, set
$$\phi(\gamma)=\frac{1}{\pi}\arg Z(\gamma)\in (\theta,\theta+1).$$
The following result gives the relation between the elements $\kappa$
and $\delta$.
\begin{prop}\label{pr:Stokes HN Reineke}
\label{prop}\hfill
\begin{enumerate}
\item For all $\gamma\in\LambdaZ_+$, there is a finite sum
\begin{equation}
\label{help}
\kappa_{\gamma}=
\sum_{n\geq 1}
\sum_{\substack
{\gamma_1+\cdots+\gamma_n=\gamma\\
\phi(\gamma_1)>\cdots>\phi(\gamma_n)}}
\delta_{\gamma_1}\cdots\delta_{\gamma_n},
\end{equation}
where the sum is over elements $\gamma_i\in \LambdaZ_+$.
\item Conversely, for $\gamma\in\LambdaZ_+$
\begin{equation}
\label{reinekeinverted}
\delta_{\gamma}=
\sum_{n\geq 1}\sum_{\stackrel
{\gamma_1+\cdots+\gamma_n=\gamma}{\phi(\gamma_1+\cdots+\gamma_i)> \phi(\gamma)}}
(-1)^{n-1}\kappa_{\gamma_1}\cdots \kappa_{\gamma_n},
\end{equation}
\end{enumerate}
\end{prop}
\begin{pf}
(i) follows from substituting \eqref{hold} into the formula of Lemma \ref{veryeasy}.
(ii) follows from Reineke's inversion of formula \eqref{help} \cite[Section 5]{Reineke}.
\end{pf}

\begin{remark}
The operation of replacing the ray $r$ by the opposite ray $-r$ exchanges
$\LambdaZ_+$ and $\LambdaZ_-$ and changes the Stokes multipliers
$(S_+,S_-)$ to $(S_-^{-1},S_+^{-1})$ thus leaving the element $\kappa$
unchanged. This gives an easy way to obtain similar expressions to \eqref
{help} and \eqref{reinekeinverted} for the case $\gamma\in \LambdaZ_-$.
\end{remark}
\begin{remark}
As pointed out in \S \ref{ss:Stokes factors}, the Stokes multipliers of $\nabla$
are determined by the Stokes factors via Lemma \ref{veryeasy}. Conversely,
it is well--known that, for $G$ reductive at least, the Stokes factors can be
recovered from the Stokes multipliers \cite[Lemma 2]{BJL}, \cite{Boalch2}.
To the best of our knowledge however, no explicit formula was known for
this procedure, even in the case of $GL_n(\IC)$. Reineke's inversion formula
\eqref{reinekeinverted} gives such a formula.
\end{remark}

\subsection{The functions $Q_n$}\label{ss:Qn}

We give below a formula for the Stokes multiplier $S_+$. The special
functions $Q_n(z_1,\ldots,z_n)$ appearing in this formula have the
property that
\[Q_n(z_1,\ldots,z_n)=0\]
unless $\pm s_n\in \Hr$, where $s_n=z_1+\cdots +z_n$. Moreover,
$Q_n(z_1,\cdots,z_n)$ is invariant under the operation of changing
$r$ to $-r$. Thus, it is enough to define $Q_n(z_1,\cdots,z_n)$
when $s_n\in \Hr$.

\begin{defn}
Set $Q_1\equiv 2\pi i$. For $n\geq 2$ and $(z_1,\ldots,z_n)\in(\IC^*)
^n$ such that $s_n\in\Hr$, define
\[Q_n(z_1,\ldots,z_n)
=2\pi i\int_{C}\frac{dt}{t-s_1}\circ \cdots \circ \frac{dt}{t-s_{n-1}}\]
where the path $C$ starts at $0$,  goes out along the ray $-r$ avoiding
any points $s_i$ by small clockwise loops, goes clockwise round a
large circle, and  finally comes back along the ray $s_n-r$, again avoiding
any points $s_j$ by small anticlockwise loops, to finish at the point $s_n$.
\end{defn}

\subsection{Formula for the Stokes multipliers}

\begin{thm}\label{one3}
If $r$ is a ray such that $\pm r$ are admissible, the components $\kappa
_\gamma$ of the Stokes multiplier $S_+$ of $\gamma$ are given by the
sum
$$\kappa_\gamma=
\sum_{n\geq 1}
\sum_{\stackrel{\alpha_i\in\PhiZ}{\alpha_1+\cdots+\alpha_n=\gamma}}
Q_n(Z(\alpha_1),\ldots, Z(\alpha_n))\,
f_{\alpha_1} f_{\alpha_2}\cdots f_{\alpha_n}$$
which is absolutely convergent over $n$.
\end{thm}
\begin{pf}
Let $i\ol{\IH}_r$ be the semi--closed half--plane $i\IH_r\sqcup -r$.
Note first that the sum in \eqref{help} may be taken over all $\gamma
_i\in\Lambda$ such that $Z(\gamma_i)\in i\ol{\IH}_r$ since, for such
$\gamma$, $\delta_\gamma=0$ unless $\gamma\in\LambdaZ_+$. 

Substituting \eqref{mainy} in \eqref{help} and using Remark \ref
{ss:Stokes map} yields, for any $\gamma\in\Lambda_+^Z$,
$$\kappa_\gamma=
\sum_{n\geq 1}
\sum_{\substack{\alpha_1,\ldots,\alpha_n\in\PhiZ\\\alpha_1+\cdots+\alpha_n=\gamma}}
\wt{Q}_n(Z(\alpha_1),\ldots,Z(\alpha_n))\,
f_{\alpha_1}\cdots f_{\alpha_n}$$
where for $z_1,\ldots,z_n\in\IC^*$ such that $s_n\in i\IH_r$, $\wt{Q}_n(z_1,\ldots,z_n)$
is defined by
\begin{equation}\label{eq:tilde Qn}
\wt{Q}_n(z_1,\ldots,z_n)
=
\sum_{\substack{
1\leq k\leq n\\[0.2 em]
0=i_0<\cdots<i_k= n\\[0.1em]
s_{i_1},s_{i_2}-s_{i_1},\ldots,s_{i_k}-s_{i_{k-1}}\in i\ol{\IH}_r\\
\phi(s_{i_1})>\phi(s_{i_2}-s_{i_1})>\cdots>\phi(s_{i_k}-s_{i_{k-1}})}}
\prod_{j=0}^{k-1} M_{i_{j+1}-i_j}(z_{i_j +1} ,\ldots,z_{i_{j+1}})
\end{equation}
with $s_j=z_1+\cdots+z_j$ and $\phi(z)=\frac{1}{\pi}\arg(z)\in (\theta,\theta+1)$.
The result now follows from the Lemma below. 
\end{pf}

\begin{lemma}\label{le:Q=Q}
For any $z_1,\ldots,z_n\in\IC^*$ such that $s_n\in i\IH_r$,
$$\wt{Q}_n(z_1,\ldots,z_n)=Q_n(z_1,\ldots,z_n)$$
\end{lemma}
\begin{pf}
To prove this claim consider the path $C$ of Definition \ref{ss:Qn} as a
piece of string and the points $s_i$ for $1\leq i\leq n-1$ as pegs. Tightening
the string will give a convex polygon with vertices some subset of the $s_i$.
Suppose $s_p$ is the first peg. Then applying Corollary \ref{corr} we can
move the string inside this first peg at the expense of adding a term
\[M_p(z_1,\cdots,z_p) Q_{n-p}(z_{p+1}, \cdots, z_n).\]
By induction, $Q_{n-p}(z_{p+1}, \cdots, z_n)$ is a sum over convex polygons
so we obtain the part of the sum on the \rhs of \eqref{eq:tilde Qn} corresponding
to $i_1=p$. Tightening the string again it catches on another peg and we repeat.
\end{pf}

\subsection{}

The following diagram summarizes the relationships between the elements $\delta,
\epsilon$ representing the Stokes factors, the element $\kappa$ representing the
Stokes multipliers, and the element $f\in\god$ defining $\nabla$.

\[\xymatrix{(f)  \ar@/_1pc/[rr]_{\stackrel{\text{\scriptsize Stokes}}{\eqref{main}}}&&(\epsilon) \ar@/_1pc/[ll]_{\stackrel{\text{\scriptsize Stokes}^{-1}}{\eqref{main2}}} \ar@/_1pc/[rr]_{\stackrel{\text{\scriptsize exp}}{\eqref{exp}}} && (\delta) \ar@/_1pc/[ll]_{\stackrel{\text{\scriptsize log}}{\eqref{log}}} \ar@/_1pc/[rr]_{\stackrel
{\text{\scriptsize clockwise multiplication}}{\eqref{help}}} &&(\kappa)\ar@/_1pc/[ll]_{\stackrel{\text{\scriptsize Reineke inversion}}{\eqref{reinekeinverted}}}
  }\]
  

\section{Proof of Proposition \ref{pr:spectral}}\label{se:Phil lemma}

\subsection{} 

Let $r$ be a ray and $\IH_r\subset\IC^*$ the open half--plane
given by
\[\IH_r=\{uv|\,u\in r, \Re(v)>0\}\]

\begin{lemma}\label{le:basic}
Let $\lambda\in\IC^*$ be a non--zero complex number. Then, the
function $e^{-\lambda/t}$ has a limit $L\in\IC$ as $t\to 0$ along the
ray $r$ if, and only if $\lambda\in \IH_r$, in which case $L=0$.
\end{lemma}
\begin{pf}
Write $\lambda=\rho e^{i\theta}$ and $t=\sigma e^{i\phi}$. Then,
$$e^{-\lambda/t}=
e^{-\frac{\rho}{\sigma}\; e^{i(\theta-\phi)}}=
e^{-\frac{\rho}{\sigma}\cos(\theta-\phi)}
e^{-i\frac{\rho}{\sigma}\sin(\theta-\phi)}$$
This has a finite limit as $\sigma\to 0$ if, and only if
$\theta\in(\phi-\frac{\pi}{2},\phi+\frac{\pi}{2})$ and, in that case,
decreases exponentially to 0.
\end{pf}

\subsection{} 

Let $U$ be a \fd vector space and $Z\in\End(U)$ a semisimple
endomorphism of $U$. Let $\sigma(Z)\subset\IC$ be the set of
eigenvalues of $Z$ and
\begin{equation}\label{eq:U eigenspaces}
U=\bigoplus_{\lambda\in\sigma(Z)}U_\lambda
\end{equation}
the corresponding decomposition of $U$ into eigenspaces of $Z$.

\begin{lemma}\label{le:linear localisation}
Let $r_1,r_2$ be two rays such that $r_1\neq -r_2$ and $u\in U$ an
element such that
$$e^{-Z/t} u\to 0
\quad\text{as}\quad
\text{$t\to 0$ in $\IH_{r_1}\cap \IH_{r_2}$}$$
Then, 
$$u\in\bigoplus_{\lambda\in\ol{\Sigma}}U_\lambda$$
where $\ol{\Sigma}\subset\IC^*$ is the closed convex sector bounded
by $r_1$ and $r_2$.
\end{lemma}
\begin{pf}
Let $u=\sum_{\lambda}u_\lambda$ be the decomposition of $u$
corresponding to \eqref{eq:U eigenspaces}. Since each $U_\lambda$
is stable under $\exp(-Z/t)$ and $e^{-Z/t}u_\lambda=e^{-\lambda/t}
u_\lambda$, we find that $u_0=0$ and that $e^{-\lambda/t}\to 0$ as
$t\to 0$ in $\IH_{r_1}\cap \IH_{r_2}$ for any $\lambda$ such that
$u_\lambda\neq 0$. Applying Lemma \ref{le:basic} to a ray $r$
contained in $\IH_{r_1}\cap\IH_{r_2}$ then shows that any such
$\lambda$ is contained in
$$\bigcap_{r\subset \IH_{r_1}\cap \IH_{r_2}}\IH_r=\ol{\Sigma}$$
\end{pf}

\subsection{} 

Let $G$ be an affine algebraic group and $\g$ its Lie algebra. Let $Z\in
\g$ be a semisimple element and decompose $\g$ as the sum $\bigoplus
_{\lambda\in\IC}\g_\lambda$ of eigenspaces for the adjoint action of $Z$.
The following result is Proposition \ref{pr:spectral} of Section
\ref{se:Irregular connections}. It was proved by Boalch \cite[lemma 6]
{Boalch2} in the case where $G$ is reductive and $Z$ is regular by using
the Bruhat decomposition of $G$.

\begin{prop}\label{pr:group localisation}
Let $r_1,r_2$ be two rays such that $r_1\neq -r_2$ and $g\in G$ an element
such that
$$e^{-Z/t}\cdot g\cdot e^{Z/t}\to 1
\quad\text{as}\quad t\to 0
\quad\text{in}\quad \IH_{r_1}\cap \IH_{r_2}$$
Then $g$ is unipotent and $X=\log(g)$ lies in 
$$\bigoplus_{\lambda\in\ol{\Sigma}}\g_\lambda$$
where $\ol{\Sigma}\subset\IC^*$ is the closed convex sector bounded
by $r_1$ and $r_2$.
\end{prop}
\begin{pf}
Embed $G$ as a closed subgroup of $GL(V)$, where $V$ is a
faithful representation. Applying Lemma \ref{le:linear localisation}
to $u=(g-1)\in\gl(V)=U$, we find that $u$ lies in the span of the
$\ad(Z)$--eigenspaces of $\gl(V)$ corresponding to eigenvalues
lying in $\ol{\Sigma}$. In particular, $u$ is a nilpotent endomorphism
of $V$. The finite sum 
\[X=\Log(g)=\sum_{n\geq 1}(-1)^{n-1}\frac{u^n}{n}\]
is then a well--defined
element of $\gl(V)$  and lies in $\g\subset\gl(V)$ because $g\in
G$. Since the $\ad(Z)$--eigenvalues of $u^n$ are contained in
the $n$--fold sum $\ol{\Sigma}+\cdots+\ol{\Sigma}\subset\ol{\Sigma}$,
the result follows.
\end{pf}


\section{Fuchsian connections and multilogarithms}
\label{polylog}

This section contains some basic results about computing parallel transport for
Fuchsian connections using iterated integrals. These results are presumably
well--known to experts but we failed to find a suitable reference.

\subsection{}

For an introduction to iterated integrals see for example \cite{Hain}. We start
by recalling their definition. Let $\omega_1,\ldots,\omega_n$ be 1--forms
defined on a domain $U\subset \C$, and $\gamma\colon[0,1]\to U$ a a path
in $U$. Let
$$\Delta=
\{(t_1,\ldots,t_n)\in [0,1]^n:0\leq t_1\leq \cdots \leq t_n\leq 1\}\subset [0,1]^n$$
be the  unit simplex. By definition,
$$\int_\gamma \omega_1\circ \cdots \circ \omega_n=
\int_\Delta f_1(t_1)\cdots f_n(t_n)\medspace dt_1\cdots dt_n$$
where $\gamma^*\omega_i=f_i(t)dt$. The following is easily checked.

\begin{lemma}\label{le:transfo rules}\hfill
\begin{enumerate}
\item Let $\ol{\gamma}(t)=\gamma(1-t)$ be the opposite path to
$\gamma$. Then,
$$\int_{\ol\gamma}\omega_1\circ \cdots \circ \omega_n=
\int_\gamma\omega_n\circ \cdots \circ \omega_1$$
\item Let $\phi:\IC\to\IC$ be a smooth map, then
$$\int_{\phi\circ\gamma}\omega_1\circ \cdots \circ \omega_n=
\int_\gamma \phi^*\omega_1\circ \cdots \circ \phi^*\omega_n$$
\end{enumerate}
\end{lemma}

There is an alternative convention obtained by using the simplex
\[\Delta^*=
\{(t_1,\ldots,t_n)\in [0,1]^n:1\geq t_1\geq \cdots \geq t_n\geq 0\}\subset [0,1]^n\]
instead of $\Delta$. We denote the resulting integral with a $*$ above the integral
sign. Thus,
$$\int^*_\gamma \omega_1\circ \cdots \circ \omega_n
=\int_{\Delta^*} f_1(t_1)\cdots f_n(t_n)\medspace dt_1\cdots dt_n$$
This convention is the more natural one for computing parallel transport and
is the one we shall use in this section. On the other hand, the convention
relying on the simplex $\Delta$ seems to be the preferred one in the study
of multilogarithms \cite{Gon,Hain}. It is of course easy to translate between
these conventions since the change of variables $t_i^*=t_{n+1-i}$ yields
\begin{equation}\label{eq:star to no star}
\int^*_\gamma \omega_1\circ \cdots \circ \omega_n=
\int_{\gamma} \omega_n\circ\cdots\circ \omega_1
\end{equation}

\subsection{}
Let $V$ be a \fd vector space and let $\poles\subset\C$ be a finite
set of points. Given a choice of residue $A_p\in \End(V)$ for each
$p\in \poles$ we can define a  meromorphic connection on  the trivial
vector bundle over $\IP^1$ with fibre $V$ by writing
\[\nablaF=d-\big.\sum_{p\in\poles} \frac{A_p}{z-p}\; dz.\]
Suppose given a smooth path $\gamma\colon [0,1]\to \C\setminus\poles$.
The parallel transport of $\nablaF$ along $\gamma$ is the invertible
linear map $\PT_\gamma\in\GL(V)$ obtained by analytically continuing
flat sections of $\nablaF$ along $\gamma$. Thus, if $\Phi$ is a fundamental
solution defined near $\gamma(0)$, then
\[\PT_{\gamma} = \Phi(\gamma(1))\cdot \Phi(\gamma(0))^{-1}.\]
where $\Phi(\gamma(1))$ is the value at $\gamma(1)$ of the analytic
continuation of $\Phi$ along $\gamma$.

Solving the differential equation for flat sections of $\nablaF$ using Picard iteration
gives the following power series expansion for such parallel transport maps (see,
\eg \cite[Lemma 2.5]{Hain}).

\begin{thm}
\label{Chen}
For any smooth path $\gamma\colon [0,1]\to \C\setminus\poles$ one has
\begin{equation*}
\PT_\gamma=
1+\Big.\sum_{\stackrel{n \geq 1}  {p_1,\ldots, p_n\in\poles}}
I_{\gamma,n}(p_1,\ldots,p_n) A_{p_1} \cdots A_{p_n},
\end{equation*}
where the sum is absolutely convergent, and the coefficients are iterated integrals
\[I_{\gamma,n}(z_1,\ldots,z_n)=
\int^*_\gamma \frac{dz}{z-z_{1}}
\circ\cdots\circ
\frac{dz}{z-z_{n}}.\]
\end{thm}

\subsection{}\label{ss:normalised limit}

We assume for the rest of the section that each of the residues $A_p$ is nilpotent.
In particular the connection $\nablaF$ is {\it non--resonant}, that is the eigenvalues
of the residues $A_p$ do not differ by positive integers. In this situation it is well--known
that for any connected and simply--connected neighbourhood $U_p$ of a pole $p
\in\poles$ there is a unique holomorphic function $H_p:U_p\rightarrow GL(V)$ with
$H_p(p)=1$ such that for any determination of the function $\log(z-p)$, the multivalued
holomorphic function
\begin{equation}\label{fun}
\Phi_p(z)=H_p(z) (z-p)^{A_p} =H_p(z)\exp( A_p\log(z-p)),
\end{equation}
is a fundamental solution of $\nablaF$. For details see, \eg \cite{Ince}. We shall refer
to $\Phi_p$ as the canonical fundamental solution of $\nablaF$ relative to a chosen
determination of $\log(z-p)$.

\begin{prop}
\label{limit}
Assume the residue $A_p$ is nilpotent and let $\Phi_p(z)$ be the canonical fundamental
solution of $\wh{\nabla}$ near $p\in\poles$. Then
\[(z-p)^{-A_p}\cdot\Phi_p(z)\to 1\text{ as } z\to p.\]
\end{prop}
\begin{pf}
Write
\[(z-p)^{-A_p} H_p(z)(z-p)^{A_p}=
H_p(z)+[(z-p)^{-A_p},H_p(z)] (z-p)^{A_p}.\]
By definition the first term tends to 1 as $z\to p$. Writing $H_p(z)=1+(z-p)J_p(z)$ with $J_p$
holomorphic at $z=p$, the second term can be rewritten as $(z-p)[(z-p)^{-A_p},J_p(z)]$. Since
$A_p$ is nilpotent and $(z-p)$ term kills all powers of $\log(z-p)$, this tends to zero as $z\to p$.
\end{pf}

\subsection{Regularised parallel transport}

Suppose now that $\gamma\colon [0,1]\to \C$ is a path such that $\gamma(0,1)\subset\C
\setminus\poles$ but which starts at a pole $p\in \poles$ and ends at a pole $q\in\poles$.
For $0<s<t<1$ let $\gamma_{[s,t]}$ denote the path in $\C\setminus\poles$ obtained by
restricting $\gamma$ to the interval $[s,t]$. It follows from Proposition \ref{limit} that the limit
\begin{equation}
\label{ho}
\PT^{\reg}_\gamma=\lim_{\stackrel{s\to 0}{t\to 1}}\bigg[ (\gamma(t)-q)^{-A_q}\cdot \PT_
{\gamma_{[s,t]}} \cdot\; (\gamma(s)-p)^{A_p}\bigg]\end{equation} is well-defined. Its value
is called the {\it regularized parallel transport} of $\nablaF$ along $\gamma$. Such limits
will be important in our computations of Stokes factors.

\begin{lemma}
\label{wellbehaved}
If $p_1\neq q$ and $p_n\neq p$ then the integral
\[I_{\gamma,n}(p_1,\ldots,p_n)=\lim_{\stackrel{s\to 0}{t\to 1}} \int^*_{\gamma_{[s,t]}} \frac{dz}{z-p_{1}}
\circ\cdots\circ
\frac{dz}{z-p_{n}}\]
is convergent.
\end{lemma}
\begin{pf}
This is proved in \cite[Section 2.9]{Gon}. 
\end{pf}

\begin{prop}\label{pr:reg transport}
Assume the residues of $\nablaF$ are nilpotent. Suppose that $P\colon U\to V$ and $Q\colon V\to W$
are linear maps such that $A_p\cdot P=0$ and $Q\cdot A_q=0$. Then,
\[Q\cdot \PT_\gamma^{\reg}\cdot P=\lim_{\stackrel{s\to 0}{t\to 1}} Q\cdot \PT_{\gamma_{[s,t]}} \cdot P,\]
and there is a series expansion
\[Q\cdot\PT_\gamma^{\reg}\cdot P=
Q\left(1+\Big.\sum_{n \geq 1}
\sum_{\stackrel{p_1,\ldots, p_n\in\poles}{p_1\neq q, p_n\neq p}}
I_{\gamma,n}(p_1,\ldots,p_n) \cdot A_{p_1} \cdots A_{p_n}\right)P\]
which is absolutely convergent in $n$.
\end{prop}

\begin{pf}
The first statement is clear because $(\gamma(s)-p)^{A_p}\cdot P=P$ and $Q\cdot
(\gamma(t)-q)^{-A_q}=Q$. To obtain the series expansion consider first fixing the
residues $A_p$ and then rescaling them by an element $\lambda\in\C$. For each
$0<s<t<1$, the function
\[Q\cdot \PT_{\gamma_{[s,t]}} \cdot P\]
is then an analytic function of $\lambda$. Theorem \ref{Chen} shows that it has Taylor
series
\[Q\cdot P+\Big.\sum_{n \geq 1} \sum_{\stackrel{p_1,\ldots, p_n\in\poles}{p_1\neq q, p_n\neq p}}
\bigg(I_{\gamma_{[s,t]},n}(p_1,\ldots,p_n) Q\cdot A_{p_1} \cdots A_{p_n}\cdot P\bigg) \lambda^n.\]
since the terms such that $p_1=q$ or $p_n=p$ are killed by $Q$ and $P$ respectively.

The following standard result of complex analysis completes the proof. Let $f_\epsilon$
be holomorphic functions on $\C$ defined for $\epsilon\in (0,1)$. Suppose that on some
closed disc $f_\epsilon\to f$ uniformly as $\epsilon\to 0$. Then $f$ is holomorphic on the
interior of the disc, and has a Taylor expansion there whose coefficients are the limits of
the Taylor coefficients of the $f_\epsilon$.
\end{pf}

\subsection{}

For later use we need one more result on regularized parallel transport maps.
Let $\alpha$ be a path starting at a point $p$ (possibly a pole of $\nablaF$)
and ending at a pole $q$, and $\beta$ a path starting from the pole $q$, both
paths otherwise avoiding the poles of $\nablaF$. The concatenation $\beta
\cdot\alpha$ can be deformed in two ways to give two paths $\gamma_+$,
$\gamma_-$ which avoid the point $q$ by a small anticlockwise (resp. clockwise)
half circle.

\begin{prop}
\label{monbig}
Assume that the residues of $\nablaF$  are nilpotent. Then
\[\PT^\reg_{\gamma_+}-\PT^\reg_{\gamma_-}=
\PT^\reg_{\beta}\cdot\, (e^{2\pi i A_q} -1)\cdot \PT^\reg_{\alpha}.\]
\end{prop}

\begin{pf}
Deforming the paths slightly we can assume that $\alpha(1-\epsilon)=\beta(\epsilon)$
for small enough $\epsilon>0$. Consider the expression
\[\PT_{\beta_{[s,1]}}\cdot \;(\beta(s)-q)^{A_q}\cdot (e^{2\pi i A_q}-1)
\cdot (\alpha(t)-q)^{-A_q} \cdot \PT_{\alpha_{[u,t]}}\cdot (\alpha(u)-p)^{A_p}.\]
Its limit as $s,u\to 0$ and $t\to 1$ is the right hand side of the stated identity.
Take $s=u=\epsilon$ and $t=1-\epsilon$. Since $z^{A_q}$ commutes with
$e^{A_q}$, the expression can be rewritten as
\[\PT_{\beta_{[\epsilon,1]}}\cdot \,(e^{2\pi i A_q} -1)\cdot
\PT_{\alpha_{[\epsilon,1-\epsilon]}}\cdot (\alpha(\epsilon)-p)^{A_p}.\]
Let $\delta$ be a small loop around $q$ starting at a point $z$. Parallel transport
of the canonical fundamental solution shows that
\[\PT_\delta=H_q(z) e^{2\pi i A_q} H_q^{-1}(z).\]
Thus the left hand side of the stated identity is
\[\PT_{\beta_{[\epsilon,1]}}\cdot\; H_q(\beta(\epsilon)) \cdot ( e^{2\pi i A_q}-1)
\cdot H_q^{-1}(\beta(\epsilon))\cdot \PT_{\alpha_{[\epsilon,1-\epsilon]}}\cdot (\alpha(\epsilon)-p)^{A_p}.\]
As $\epsilon\to 0$, $\beta(\epsilon)\to q$ and $H_q(\beta(\epsilon))\to 1$. This
gives the result.
\end{pf}

\subsection{}

Let $(z_1,\ldots,z_n)\in\C^n$, $z\in\C$ and assume that $z=z_i$ for a unique
$i$. Assume that $\alpha$ is a path ending at $z$ and $\beta$ is a path starting
at $z$, both paths otherwise avoiding the points $z_i$. Let $\gamma_+$ (resp.
$\gamma_-$ be the paths obtained by deforming the concatenation $\beta\cdot
\alpha$ by avoiding the point $z_i$ by a small small anticlockwise (resp. clockwise)
half circle.

\begin{cor}\label{corr}
The following holds
\begin{multline*}
I_{\gamma_+,n}(z_1,\ldots,z_n)-I_{\gamma_-,n}(z_1,\ldots,z_n)\\
=2\pi i \cdot I_{\alpha,i-1}(z_1,\ldots,z_{i-1})\cdot I_{\beta,n-i}(z_{i+1},\ldots,z_n).
\end{multline*}
\end{cor}
\begin{pf}
This is proved in \cite[Cor. 2.6]{Gon}. Alternatively it can easily be obtained by
equating coefficients in both sides of Proposition \ref{monbig}.
\end{pf}


\section{Fourier--Laplace transform}
\label{se:Fourier-Laplace}

In this section, we prove Theorem \ref{jurk} for the general linear groups
by using the \FL transform of the connection $\nabla$.

\subsection{} \label{ss:Z F}

Let $V$ be a complex, \fd vector space, $\V$ the holomorphically trivial
vector bundle on $\IP^1$ with fibre $V$ and $\nabla^V$ the meromorphic
connection on $\V$ given by
\begin{equation}\label{eq:nabla V}
\nabla^V=d-\left(\frac{Z}{t^2}+\frac{F}{t}\right)dt
\end{equation}
where $Z,F\in\gl(V)$. The assumptions $(Z)$--$(f)$ of Section \ref{ss:assumptions}
translate into the following ones:
\begin{enumerate}
\item[(Z)] $Z$ is diagonalisable. We denote the roots of the minimal polynomial
of $Z$ by $z_1,\ldots,z_m$, the corresponding eigenspaces by $V_1,\ldots,
V_m$ and let $P_1,\ldots,P_m$ be the projections corresponding to the
decomposition
\begin{equation}\label{eq:eigenspaces}
V=V_1\oplus\cdots\oplus V_m
\end{equation}
\item[(F)] The diagonal blocks of $F$ with respect to the decomposition
\eqref{eq:eigenspaces} are zero, that is $P_iFP_i=0$ for any $i$.
\end{enumerate}
According to Definition \ref{de:Stokes rays}, the Stokes rays of $\nabla
^V$ are the rays $\ray{(z_i-z_j)}$, $1\leq i\neq j\leq m$.

\subsection{}\label{ss:FL}

Let $\wh{\V}$ be another copy of the trivial vector bundle on $\IP^1$ with
fibre $V$ and consider the Fuchsian connection $\nablah^V$ on $\wh{\V}
$ with poles at the points $z_1,\ldots,z_m$ given by
$$\nablah^V=
d-\sum_{i=1}^m\frac{P_i F}{z-z_i}\;dz$$
The connection $\nablah^V$ is of the form considered in Section \ref{polylog}.
Moreover, since
$$(P_i F)^2=(P_i F P_i)F=0$$
by assumption $(F)$, the residues $A_i=P_i F$ are nilpotent. In particular,
$\nablah^V$ is non--resonant.

\subsection{}

Fix a pole $z_i$ and let $Q_i:V_i\hookrightarrow V$ be the inclusion. Let
$U_i$ be a connected and simply--connected neighborhood of $z_i$ in
$\IP^1\setminus\{z_1,\ldots,\wh{z_i},\ldots,z_m\}$.

\begin{lemma}
\label{drive}
There is a unique horizontal section $\phi^{(z_i)}$ of $\nablah^V$ defined
on $U_i$ and taking values in $\Hom_{\IC}(V_i,V)$ which is regular at $z_i$
and such that $\phi^{(z_i)}(z_i)$ is the inclusion $Q_i:V_i\into V$.
\end{lemma}
\begin{pf}
Let $\Phi_i:U_i\to GL(V)$ be the canonical fundamental solution of $\nablah
^V$ at $z_i$ (see \S \ref{ss:normalised limit}). Since $P_iFQ_i=0$ by assumption
$(F)$,
$$\phi^{(z_i)}(z):=\Phi_i(z)\cdot Q_i=
H_i(z)(z-z_i)^{P_iF}\cdot Q_i=H_i(z)\cdot Q_i$$
gives the required solution. Uniqueness is straightforward.
\end{pf}

\subsection{}

Fix a pole $z_i$, an admissible ray $r=\ray{e^{i\varphi}}$, and set
\begin{equation}\label{eq:Y_ell}
Y_r^{(z_i)}(t)=\frac{1}{t}\int_{z_i+r}\phi^{(z_i)}(z)e^{-z/t}dz
\end{equation}

\begin{prop}\label{pr:Fourier}\hfill
\begin{enumerate}
\item The integral \eqref{eq:Y_ell} is convergent for any $t$ in the
half--plane $\IH_r$.
\item The corresponding function
$Y^{(z_i)}_r\colon \IH_r\to \Hom_{\C}(V_i,V)$
is holomorphic and satisfies
\begin{equation}\label{eq:irregular}
\frac{dY_r^{(z_i)}}{dt}=\left(\frac{Z}{t^2}+\frac{F}{t}\right)Y_r^{(z_i)}
\end{equation}
\item $Y^{(z_i)}_r\cdot e^{z_i/t}$ tends to the inclusion $Q_i\colon
V_i\into V$ as $t\to 0$ in $\IH_r$.
\end{enumerate}
\end{prop}
\begin{pf}
We drop the superscript $(z_i)$ and the subscripts $r$ and $z_i+
r$ from the notation and use primes for derivatives.
(i) follows from the fact that since $\nablah^V$ has regular singularities,
$\phi^{(z_i)}$ grows at most polynomially as $z\rightarrow\infty$, while
$$|e^{-z/t}|=|e^{-z_i/t}|\,e^{-\frac{|z-z_i|}{|t|}\cos(\varphi-\arg(t))}$$
decreases exponentially as $z\rightarrow\infty$ along $z_i+\ray{e^{i\varphi}}$,
provided $\arg(t)\in(\varphi-\frac{\pi}{2},\varphi+\frac{\pi}{2})$.

(ii) Differentiating the defining integral for
$Y(t)$ gives
\[Y'(t)=-\frac{1}{t} Y(t)+\frac{1}{t^3}\int \phi(z)z e^{-z/t} dz.\]
Integrating the second term by parts gives
\[Y'(t)=-\frac{1}{t} Y(t)+
\frac{1}{t^2}Q_i z_ie^{-z_i/t}+
\frac{1}{t^2}\int\frac{d}{dz}\big(z\phi(z)\big) e^{-z/t} dz.\]
Expanding the derivative, two terms cancel, giving
\[Y'(t)=\frac{1}{t^2}Q_i z_ie^{-z_i/t}+
\frac{1}{t^2}\int
\bigg(\sum_j\frac{P_j F}{z-z_j}\,\phi(z)\bigg) z e^{-z/t} dz.\]
Taking the finite sum outside the integral and writing $z=(z-z_j)+z_j$
gives a sum of two terms. The first is
\[\frac{1}{t^2}\sum_j P_j F\,\int \phi(z)e^{-z/t}dz=
\frac{1}{t}\bigg(\sum_j P_j \bigg)F\,Y(t)=
\frac{F}{t}Y(t).\]
The second is
\begin{equation}\label{eq:second}
\frac{1}{t^2}\left(Q_i z_ie^{-z_i/t}+
\sum_j  z_j P_j F\,
\int\frac{\phi(z)e^{-z/t}}{z-z_j}\;dz\right).
\end{equation}
To simplify this expression, we integrate by parts the defining expression
for $Y(t)$ to get
\begin{equation*}
\begin{split}
Y(t)
&=Q_i e^{-z_i/t}+\int \phi'(z) e^{-z/t} dz\\
&=Q_i e^{-z_i/t}+
\sum_j P_j F\,\int \frac{\phi(z)e^{-z/t}}{z-z_j}\;dz
\end{split}
\end{equation*}
which, upon being multiplied by $P_j$ yields
$$P_j Y(t)=\delta_{ij}\,Q_i e^{-z_i/t}+
P_j F\,\int \frac{\phi(z)e^{-z/t}}{z-z_j}\;dz$$
Substituting this into \eqref{eq:second} shows that the latter is equal to
$$\frac{1}{t^2}\sum_j z_j\,P_j Y(t)=\frac{Z}{t^2} Y(t)$$
and therefore that $Y$ satisfies \eqref{eq:irregular} as claimed.

(iii) The limiting behaviour of $Y(t)$ as $t\to 0$ in $\IH_r$ follows at once
from the fact that $Y$ has an asymptotic expansion in $\IH_r$ with constant
term $1$. This in turn is a consequence of the Lemma below.
\end{pf}

\begin{lemma}[Watson \cite{Hardy}]
Set
\[Y(t)=\int_0^\infty  \phi(z) e^{-z/t} dz,\]
and suppose that $\phi(z)$ is analytic at $z=0$ with Taylor series
\[\phi(z)=\sum_{n\geq 0} \frac{a_n}{n!}z^n.\]
Then $Y(t)$ has the asymptotic expansion
\[Y(t)\sim \sum_{n=0}^\infty {a_n }{t^{n+1}}.\]
as $t\to 0$ in the half--plane $\Re(t)>0$.\qed
\end{lemma}

\subsection{Proof of Theorem \ref{jurk} for $\mathbf{G=GL(V)}$}
\label{ss:existence and uniqueness}

Let the functions $Y^{(z_i)}_r\colon \IH_r\to\Hom(V_i,V)$ be given
by Proposition \ref{pr:Fourier} and define a map $Y^V_r\colon \IH_r
\to\End(V)$ by
$$Y^V_r=\sum_i Y^{(z_i)}_r P_i$$
By Proposition \ref{pr:Fourier}, $Y^V_r(t)$ is a fundamental solution
of $\nabla^V$ and $Y^V(t)\cdot e^{Z/t}$ tends to the identity as $t
\rightarrow 0$ in $\IH_r$.

To prove uniqueness, let $Y_1,Y_2:\IH_r\to GL(V)$ be two holomorphic
functions satisfying \eqref{eq:diff equ}--\eqref{eq:asym}. Then $g=Y_2
^{-1}\cdot Y_1$ is a locally constant and therefore constant $GL(V)
$--valued function on $\IH_r$ such that $e^{-Z/t}\cdot g\cdot e^{Z/t}$
tends to $1$ as $t\to 0$ in $\IH_r$. Applying Proposition \ref{pr:group localisation},
we see that $g=\exp(X)$ where $X\in\gl(V)$ only has components along
the $\ad(Z)$--eigenspaces corresponding to eigenvalues $\lambda\in r$.
Since $r$ is an admissible ray of $\nabla^V$, $r$ contains no such
eigenvalues and $X=0$.\qed


\section{Computation of the Stokes factors for $GL(V)$}\label{se:Stokes factors}

We now explicitly compute the Stokes factors of the connection
$\nabla^V$ in terms of multilogarithms. Retain the notation of
Section \ref{se:Fourier-Laplace}.

\subsection{}

Fix a pole $z_i$ and let $\ell$ be a Stokes ray of $\nabla^V$, so that
$z_i+\ell$ contains some of the poles of $\nablah^V$. List these in
order of increasing distance \wrt $z_i$ as $z_{i_j}$, $j=1,\ldots,p$ and
set $z_{i_0}=z_i$. Let $r_\pm$ be small anticlockwise and clockwise
rotations of $\ell$ respectively such that the closed, convex sectors
$\ol{\Sigma}_\pm$ determined by $z_i+r_\pm$ and $z_i+\ell$ only
contain the poles $z_{i_j}$, $j=0,\ldots,p$.

\begin{prop}\label{pr:Stokes factors}
The following holds on $\IH_{r_+}\cap\IH_{r_-}$
\begin{equation}
Y^{(z_i)}_{r_-}=Y^{(z_i)}_{r_+}
+2\pi i
\sum_{j=1}^p Y^{(z_{i_j})}_{r_+}\circ P_{i_j}F\circ \PT^\reg_{C_{i_j i}}\circ Q_i
\end{equation}
where $C_{i_j i}$ is a small perturbation of the oriented line segment
$[z_i,z_{i_j}]$ which avoids the poles $z_{i_k}$, with $1\leq k\leq j-1$
by using anticlockwise arcs of circle around them.
\end{prop}

\begin{pf}
For any $0\leq j\leq p$, let $\gamma_j$ be a small perturbation of the ray
$z_i+\ell$ which avoids the poles $z_{i_k}$, with $1\leq k\leq j$, by going
into $\ol{\Sigma}_-$ and the poles $z_{i_k}$, with $j+1\leq k\leq p$ by going
into $\ol{\Sigma}_+$. By Cauchy's theorem, and the fact that $\phi^{(z_i)}
(z)e^{-z/t}$ decays exponentially as $z$ goes to infinity along $\ell$,
$$Y^{(z_i)}_{r_+}(t)=\frac{1}{t}\int_{\gamma_0}\phi^{(z_i)}(z)e^{-z/t}dz
\aand
Y^{(z_i)}_{r_-}(t)=\frac{1}{t}\int_{\gamma_p}\phi^{(z_i)}(z)e^{-z/t}dz$$
We wish to compute the difference $\int_{\gamma_j}\phi^{(z_i)}(z)e^{-z/t}
dz-\int_{\gamma_{j-1}}\phi^{(z_i)}(z)e^{-z/t}dz$ for any $p\geq j\geq 1$.
Let $\epsilon>0$ be small enough and write
$$\gamma_j=
\gamma_{j,\epsilon}^+\circ D_{j,\epsilon}^-\circ\gamma_{j,\epsilon}^-
\aand
\gamma_{j-1}=
\gamma_{j,\epsilon}^+\circ D_{j,\epsilon}^+\circ\gamma_{j,\epsilon}^-$$
where:
\begin{itemize}
\item $\gamma_{j,\epsilon}^-$ is the perturbation of the straight line
segment from $z_i$ to $z_{i_j}-\epsilon e^{i\varphi}$, where $\ell=\IR_{>0}
\cdot e^{i\varphi}$, which avoids the poles $z_{i_k}$, with $1\leq k\leq j-1$,
by going into $\ol{\Sigma}_-$.
\item $D_{j,\epsilon}^\pm$ are the arcs of circle of radius $\epsilon$ centred
at $z_{i_j}$ joining $z_{i_j}-\epsilon e^{i\varphi}$ to $z_{i_j}+\epsilon
e^{i\varphi}$ in $\ol{\Sigma}_\pm$.
\item $\gamma_{j,\epsilon}^+$ is the perturbation of the line $z_{i_j}+\IR_
{>\epsilon}\cdot e^{i\varphi}$ which avoids the poles $z_{i_k}$, with $j+1
\leq k \leq p$, by going into $\ol{\Sigma}_+$.
\end{itemize}
This yields
\begin{multline*}
\int_{\gamma_j}\phi^{(z_i)}(z)e^{-z/t}dz-
\int_{\gamma_{j-1}}\phi^{(z_i)}(z)e^{-z/t}dz\\
=
\int_{D_{j,\epsilon}^-}\phi^{(z_i)}(z)e^{-z/t}dz-
\int_{D_{j,\epsilon}^+}\phi^{(z_i)}(z)e^{-z/t}dz+
\int_{\gamma_{j,\epsilon}^+}\left(\phi^{(z_i)}_--\phi^{(z_i)}_+\right)(z)e^{-z/t}dz
\end{multline*}
where $\phi^{(z_i)}_\mp(z)$ is the analytic continuation of $\phi^{(z_i)}$ along
the path
$$C_{j,\epsilon}^\mp(z)=
\gamma_{j,\epsilon}^+(z)\circ D_{j,\epsilon}^\mp\circ\gamma_{j,\epsilon}^-$$
and $\gamma_{j,\epsilon}^+(z)$ is the portion of $\gamma_{j,\epsilon}^+$ joining
$z_{i_j}+\epsilon e^{i\varphi}$ to $z$. Since $C_{j,\epsilon}^\mp(z)$ differ by a
small loop around $z_{i_j}$, Proposition \ref{monbig} yields
$$\int_{\gamma_{j,\epsilon}^+}\left(\phi^{(z_i)}_--\phi^{(z_i)}_+\right)(z)e^{-z/t}dz=
2\pi i\int_{\gamma_{j,\epsilon}^+}\phi^{(z_{i_j})}e^{-z/t}dz\circ
P_{i_j} F\circ \PT_{C_{i_j i}}\circ Q_i$$
By Cauchy's theorem again,
$$\lim_{\epsilon\to 0}\int_{\gamma_{j,\epsilon}^+}\phi^{(z_{i_j})}e^{-z/t}dz=
\int_{z_{i_j}+r_+}\phi^{(z_{i_j})}e^{-z/t}dz=tY^{(z_{i_j})}_{r_+}(t)$$
To conclude, it suffices to show that the integrals $\int_{D_{j,\epsilon}^\pm}
\phi^{(z_i)}(z)e^{-z/t}dz$ tend to zero as $\epsilon\to 0$. Let $\Phi_{i_j}(z)=
H_{i_j}(z)(z-z_{i_j})^{P_{i_j}F}$ be the canonical fundamental solution of
$\nablahv$ at $z_{i_j}$, so that $\phi^{(z_i)}=\Phi_{i_j}\cdot C_i$ for some
$C_i\in\Hom_\IC(V_i,V)$. Then, if $\|\cdot\|$ is an algebra norm on $\End(V)$,
\begin{multline*}
\left\|\int_{D_{j,\epsilon}^\pm}\phi^{(z_i)}(z)e^{-z/t}dz\right\|=
\left\|\int_{D_{j,\epsilon}^\pm}H_{i_j}(z)(z-z_{i_j})^{P_{i_j}F}C_ie^{-z/t}dz\right\|\\
\leq
\pi\epsilon(1+(|\ln\epsilon|+\pi)\|P_{i_j}F\|)M
\end{multline*}
where $M=\|C_i\|\cdot\max_{|z-z_{i_j}|\leq\epsilon}\|H_{i_j}(z)\|e^{-z/t}<\infty$.
\end{pf}

\subsection{}

Let $\ell$ be the Stokes ray of $\nabla^V$ and $S^V_\ell\in\GL(V)$ the
corresponding Stokes factor, so that, on $\IH_{r_-}\cap \IH_{r_+}$
\begin{equation}\label{eq:V Stokes}
Y^V_{r_-}=Y^V_{r_+}\cdot S^V_\ell
\end{equation}
where $r_{\pm}$ are small anticlockwise and clockwise perturbations of
$\ell$. Define a partial order on the set of poles of $\nabla^V$ by
\begin{equation}\label{eq:ray order}
z_j>_\ell z_i\quad\text{if}\quad z_j\in z_i+\ell
\end{equation}
The following result gives a formula for the blocks 
of $S_\ell^V$ corresponding to the decomposition \eqref{eq:eigenspaces}
of $V$ into eigenspaces of $Z$.

\begin{thm}\label{th:V Stokes}
\begin{equation}\label{eq:compact}
P_j\circ S^V_\ell\circ P_i=
\left\{\begin{array}{cll}
0	&\text{if $z_j\ngeq_\ell z_i$}\\
P_i	&\text{if $z_j=z_i$}\\
2\pi i\cdot P_j F\circ \PT^\reg_{C_{ji}}\circ P_i&\text{if $z_j>_\ell z_i$}
\end{array}\right.
\end{equation}
where $C_{ji}$ is a small perturbation of the oriented line segment $[z_i,
z_j]$ which avoids the poles $z_k\in (z_i,z_j)$ by using anticlockwise arcs
of circle around them.

Thus, for $z_j>_\ell z_i$,
$$P_j\circ S^V_\ell\circ P_i=
2\pi i\cdot P_jF\left(1+
\sum_{n\geq 1}\sum_{\substack{k_1,\ldots,k_n:\\k_1\neq j,k_n\neq i}}
\int^*_{C_{ji}}\frac{dz}{z-z_{k_1}}\circ\cdots\circ\frac{dz}{z-z_{k_n}}
P_{k_1}F\cdots P_{k_n}F\right)
P_i$$
\end{thm}
\begin{pf}
The first statement follows from \eqref{eq:V Stokes}, the fact that
$Y^V_{r_\pm}\cdot P_i=Y^{(z_i)}_{r_\pm}\cdot P_i$ and Proposition
\ref{pr:Stokes factors}. The second from Proposition \ref{pr:reg transport}.
\end{pf}

\remark Theorem \ref{th:V Stokes} and its proof extend the computation of
\BJL \cite{BJL2} to the case when $Z$ has repeated eigenvalues.


\section{Tannaka duality}\label{se:Tannaka}

In this section, we prove Theorems \ref{jurk} and \ref{one2} for
an arbitrary algebraic group $G$ by relying on the results of
Sections \ref{se:Fourier-Laplace} and \ref{se:Stokes factors}
and using Tannaka duality.

\subsection{}

Let $\rho:G\to\GL(V)$ be a \fd representation of $G$ and $\V
=P\times_G V$ the holomorphically trivial vector bundle over
$\IP^1$ with fibre $V$. The connection \eqref{nab} induces a
meromorphic connection $\nabla^V$ on $\V$ given by
\begin{equation}\label{eq:induced}
\nabla^V=d-\left(\frac{\rho(Z)}{t^2}+\frac{\rho(f)}{t}\right)dt
\end{equation}

\begin{lemma}
The connection $\nabla^V$ satisfies the assumptions of Section
\ref{ss:Z F}, that is
\begin{itemize}
\item
$\rho(Z)\in\gl(V)$ is semisimple.
\item
The diagonal blocks of $\rho(f)$ \wrt the eigenspace
decomposition of $V$ under $\rho(Z)$ is zero.
\end{itemize}
\end{lemma}
\begin{pf}
This holds because $\g=\g^Z\oplus[Z,\g]$ and $\rho:\g\to\gl(V)$
is equivariant with respect to the adjoint action of $Z$.
\end{pf}

\remark\label{rk:V Stokes} The Stokes rays of $\nabla$ and $\nabla
^V$ are related, though not in an entirely straightforward way. If $\ell$
is a Stokes ray of $\nabla$, and the restriction of $\rho$ to the subalgebra
$\bigoplus_{\zeta\in\ell}\g_\zeta\subset\g$ is not zero, then $\ell$ is also
a Stokes ray of $\nabla^V$. This is the case if $\rho$ is faithful as a
representation of $\g$ for example, and therefore if $G$ is semisimple.
In general however, simple examples show that a Stokes ray of $\nabla$
need not be one of $\nabla^V$ and that, conversely, a Stokes ray of
$\nabla^V$ need not be one of $\nabla$.

\subsection{}

Let $\RRep(G)$ be the category of \fd representations of $G$. We establish
below the naturality of the canonical fundamental solutions of the connections
$\nabla^V$ \wrt tensor products and homomorphisms in $\RRep(G)$.

The union $\SL$ of the sets of Stokes rays of the connections $\nabla^V$,
$V\in\RRep(G)$ is at most countable. Fix $r\notin\SL$ and let $\IH_r\subset
\IC^*$ be the corresponding half--plane \eqref{eq:halfplane}.

\begin{prop}\label{pr:naturality}
Let $\{Y_r^V\}_{V\in\RRep(G)}$ be a family of holomorphic functions $Y_r^V:
\IH_r\to GL(V)$ such that
\begin{gather*}
\frac{dY_r^V}{dt}=\left(\frac{\rho(Z)}{t^2}+\frac{\rho(f)}{t}\right)Y_r^V\\[1.1ex]
Y_r^V\cdot e^{\rho(Z)/t}\to 1\quad\text{as}\quad\text{$t\to 0$ in $\IH_r$}
\label{eq:cond 2}
\end{gather*}
Then, the following holds for any $V_1,V_2\in\RRep(G)$ and $T\in\Hom_G
(V_1,V_2)$,
\begin{gather}
T\,Y^{V_1}_r = Y^{V_2}_r\,T \label{eq:hom}\\[1.1ex]
Y^{V_1\otimes V_2}_r = Y^{V_1}_r\otimes Y^{V_2}_r \label{eq:tensor}
\end{gather}
\end{prop}
\begin{pf}
\eqref{eq:hom} follows from the uniqueness part of Theorem \ref{jurk}
for the group $GL(V_1\otimes V_2)$ (see \S \ref{ss:existence and uniqueness})
since both sides are fundamental solutions of $\nabla^{V_1\otimes V_2}
=\nabla^{V_1}\otimes 1+1\otimes \nabla^{V_2}$ having the required
asymptotic properties on $\IH_r$.
\eqref{eq:tensor} follows in a similar manner. Namely, consider the element
$$C=(Y_r^{V_2})^{-1}\cdot T Y_r^{V_1}\in\Hom_{\IC}(V_1,V_2)$$
Condition \eqref{eq:cond 2} implies that $e^{-
\rho_{V_2}(Z)/t}\cdot(C-T)\cdot e^{\rho_{V_1}(Z)/t}$ tends to $0$ as
$t\to 0$ in $\IH_r$. Applying Lemma \ref{le:linear localisation} to
$u=C-T\in\Hom_{\IC}(V_1,V_2)=U$, we see that the only non--trivial
components of $C-T$ along the eigenspace decomposition of the
$G$--module $\Hom_\IC(V_1,V_2)$ under $Z$ correspond to eigenvalues
lying in $r$. Since $r$ is not a Stokes ray of $\nabla^{\End(V_1,V_2)}$
however it follows that $C-T=0$.
\end{pf}

\subsection{Proof of Theorem \ref{jurk}}

Assume first that the ray $r$ is admissible for all connections $\nabla^V$,
$V\in\RRep(G)$, that is that $r\notin\SL$. For any $V\in\RRep(G)$, let $Y
_r^V:\IH_r\to\GL(V)$ be the corresponding canonical fundamental solution
of $\nabla^V$. By \eqref{eq:hom}, the collection $\{Y_r^V(t)\}$ defines a
function $Y_r$ on $\IH_r$ with values in $\wh{U\g}$. By \eqref{eq:tensor},
$Y_r$ takes values in $G\subset\wh{U\g}$ since, by Tannaka duality $G$
is the set of grouplike elements of $\wh{U\g}$ \cite{SR}. Since $Y_r^V$
satisfies the properties \eqref{eq:diff equ}--\eqref{eq:asym} for any
representation $V$ of $G$, $Y^Y_r$ is a holomorphic $G$--valued
function which satisfies these same properties.

Assume now that $r$ is admissible for $\nabla$ and let $r_+,r_-$ be small
anticlockwise and clockwise perturbations of $r$ such that $r_\pm\notin\SL$
and the closed convex sector $\ol{\Sigma}\subset \IC^*$ bounded by $r_\pm$
does not contain any Stokes ray of $\nabla$. The element $g\in G$ defined by $Y_{r_-}=Y_{r_+}g$
on $\IH_{r_-}\cap\IH_{r_+}$ is such that $e^{-Z/t}ge^{Z/t}\to 1$ as $t\to 0$ in $\IH_
{r_-}\cap\IH_{r_+}$. By Proposition \ref{pr:spectral}, $g=\exp(X)$ where $X$
lies in the span of the eigenspaces of $\ad(Z)$ corresponding to eigenvalues
contained in $\ol{\Sigma}$. Since there are none, $X=0$ so that $Y_{r_{\pm}}$
patch to a fundamental solution of $\nabla$ having the required asymptotic
property on $\IH_{r_-}\cup\IH_{r_+}\supset\IH_r$.\qed

\subsection{Proof of Theorem \ref{one2}}\label{ss:factors for G}

Let $\rho:G\to\GL(V)$ be a \fd representation of $G$. We begin by
reworking the formula for the Stokes factors of the linear connection
$\nabla^V$ obtained in Theorem \ref{th:V Stokes}.

Let $z_1,\ldots,z_m$ be the roots of the minimal polynomial of $\rho
(Z)$ and $P_1,\ldots,P_m\in\End(V)$ the corresponding eigenprojections.
Let $\ell$ be a Stokes ray of $\nabla^V$ and $z_i,z_j$ two eigenvalues
of $\rho(Z)$ such that $z_j\in z_i+\ell$. Since $f=\sum_{\alpha\in\PhiZ}
f_\alpha$, the product of operators $P_jfP_{k_1}f\cdots P_{k_n}fP_i$
appearing in Theorem \ref{th:V Stokes} is equal to
$$\sum_{\alpha_0,\ldots,\alpha_n\in\PhiZ}
P_j f_{\alpha_0} P_{k_1}f_{\alpha_1}
\cdots P_{k_n}f_{\alpha_n} P_i$$
where we abusively denote $\rho(f)$ by $f$.

Given that $P_{k_m}f_{\alpha_m}P_{k_{m+1}}=0$ unless $z_{k_m}
=z_{k_{m+1}}+Z(\alpha_m)$, the sum over $k_1,\ldots,k_n$ becomes
one over the roots $\alpha_0,\ldots,\alpha_n$ with $z_{k_m}=z_j-
Z(\alpha_0+\cdots+\alpha_{m-1})$ and the constraint $Z(\alpha_0+
\cdots+\alpha_m)=z_j-z_i$. It follows that $(2\pi i)^{-1}P_jS_\ell^VP_i$
is equal to
\[P_jfP_i+
\sum_{\substack{\alpha_0,\ldots,\alpha_n\in\PhiZ:\\Z(\alpha_0+\cdots+\alpha_n)=z_j-z_i}}
\int_{C_{ji}}^*
\frac{dz}{z-z_j+s_0}\circ\cdots\circ\frac{dz}{z-z_j+s_{n-1}}\cdot
P_j f_{\alpha_0} f_{\alpha_1}\cdots f_{\alpha_n} P_i\]
where $s_m=Z(\alpha_0+\cdots+\alpha_m)$, so that $s_n=z_j-z_i$.
By Lemma \ref{le:transfo rules}, the change of variable $z\to z_j-z$
in the iterated integral yields
$$\int_{C_{0,s_n}}^*
\frac{dz}{z-s_0}\circ\cdots\circ\frac{dz}{z-s_{n-1}}
=
\int_{\ol{C_{0,s_n}}}
\frac{dz}{z-s_0}\circ\cdots\circ\frac{dz}{z-s_{n-1}}$$
where $C_{0,s_n}$ is a small perturbation of $[s_n,0]$ which avoids the poles
$s_k$, $1\leq k\leq n-1$ such that $s_k\in(s_n,0)$ by using anticlockwise arcs
of circles around them, $\ol{C_{0,s_n}}(t)=C_{0,s_n}(1-t)$ is $C_{0,s_n}$ with
the opposite orientation, and the above equality follows from (i) of Lemma \ref
{le:transfo rules} and \eqref{eq:star to no star}. Since $\ol{C_{0,s_n}}$ is a
perturbation of $[0,s_n]$ which avoids the $s_k$ such that
$s_k\in(s_n,0)$ by {\it clockwise} arcs of circles around them, it follows that
\[P_jS_\ell^VP_i=
2\pi i
\sum_{n\geq 1}
\sum_{\substack{\alpha_1,\ldots,\alpha_n\in\PhiZ:\\Z(\alpha_1+\cdots+\alpha_n)=z_j-z_i}}
M_n(Z(\alpha_1),\ldots,Z(\alpha_n))
P_jf_{\alpha_1}\cdots f_{\alpha_n}P_i\]
where the $M_n$ are the functions defined in Section \ref{ss:formula for factors}.
Since $S_\ell^V$ is equal to $1+\sum_{i,j:z_j-z_i\in\ell}P_jS_\ell^VP_i$, 
this yields
\begin{equation}\label{eq:almost there}
S_\ell^V=
1+2\pi i
\sum_{n\geq 1}
\sum_{\substack{\alpha_1,\ldots,\alpha_n\in\PhiZ:\\Z(\alpha_1+\cdots+\alpha_n)\in\ell}}
M_n(Z(\alpha_1),\ldots,Z(\alpha_n))
f_{\alpha_1}\cdots f_{\alpha_n}
\end{equation}

Let now $\ell$ be a Stokes ray of $\nabla$. We wish to show that $\rho(S_
\ell)$ is given by the \rhs of \eqref{eq:almost there}. If $\ell$ is also a Stokes
ray of $\nabla^V$, this follows from the fact that $\rho(S_\ell)=S_\ell^V$. If,
on the other hand, $\ell$ is not a Stokes ray of $\nabla^V$, then $\rho(S_\ell)=1$
which is also the value returned by the \rhs of \eqref{eq:almost there}. Indeed,
in this case $z_j-z_i\notin\ell$ for any $i,j$, so that $P_jf_{\alpha_1}\cdots
f_{\alpha_n}P_i=0$ whenever $Z(\alpha_1+\cdots+\alpha_n)\in\ell$.
\qed


\section{Inversion of non--commutative power series}
\label{comb}

In this section, we compute the Taylor series of the inverse of the Stokes map,
thereby proving Theorem \ref{th:Stokes inverse}. We shall do so by working out
a non--commutative analogue of the compositional inversion of a formal power
series.

\subsection{} 

We shall need some notation on trees. A tree $T$ is a finite, connected and
simply--connected graph.  We denote the set of edges of $T$  by $E(T)$ and
the set of vertices by $V(T)$. A plane tree is a tree $T$ together with a cyclic
ordering of the incident edges at each vertex. A tree has both internal and
external edges; a rooted plane tree is a plane tree with a distinguished
external edge called the root; the other external edges are then called the
leaves.  We draw plane trees in such a way that the cyclic ordering of the
edges incident at a given vertex is the natural clockwise ordering induced
by the embedding in the plane. For example
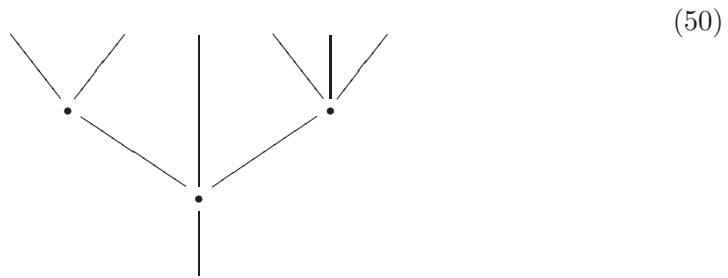
\begin{equation}\label{eq:tree}
\xymatrix@C=1.5em
{\ar@{-}[dr] &&   \ar@{-}[dl] &  \ar@{-}[dd] &  \ar@{-}[dr] & \ar@{-}[d] & \ar@{-}[dl] \\
                   &\blob \ar@{-}[drr] &&&&\blob \ar@{-}[dll] \\ 
                   &&&\blob \ar@{-}[d] \\
                   &&& \
}
\end{equation}
Note that the incoming edges at a vertex $v$ have a canonical ordering.
Similarly the leaves of $T$ have a canonical ordering.

\subsection{}\label{ss:NC series} 

Let $U$ and $V$  be  complex vector spaces. By a {\it non--commutative
(NC) power series} $\phi\colon U\to V$ we mean a sequence of linear maps
\[\phi_n\colon U^{\tensor n}\to V,\quad {n\geq 1}.\]
Two such power series $\phi\colon U\to V$ and $\psi\colon V\to W$ can
be composed to give a power series $\psi\circ\phi\colon U\to W$ by using
the following rule
\begin{multline}
\label{comp}
(\psi\circ\phi)_n(u_1,\ldots,u_n)= \\
\sum_{k=1}^n\thickspace\thickspace
\sum_{{0=i_0<\cdots <i_k=n}}
\psi_k \big(\phi_{i_1-i_0}(u_{i_0+1},\ldots, u_{i_1}), \cdots,
\phi_{i_k-i_{k-1}}(u_{i_{k-1} +1},\ldots, u_{i_k})\big).\end{multline}
This sum is best visualized as a sum over plane rooted trees of height two with
the tensors  $\phi$ and $\psi$  labelling the  vertices, and the inputs $u_1,\ldots,
u_n$ labelling the leaves. For example, in the sum for $(\psi\circ \phi)_6$, the term
\[\psi_3(\phi_2(u_1,u_2),\phi_1(u_3),\phi_3(u_4,u_5,u_6))\]
corresponds to the tree
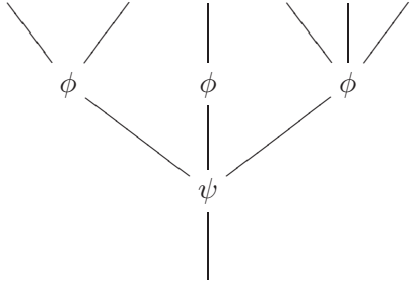
\[
\xymatrix@C=1.5em{  \ar@{-}[dr] &&   \ar@{-}[dl] &  \ar@{-}[d] &  \ar@{-}[dr] & \ar@{-}[d] & \ar@{-}[dl] \\
                   &\phi \ar@{-}[drr] && \phi \ar@{-}[d] &&\phi \ar@{-}[dll] \\
                   &&&\psi \ar@{-}[d] \\
                   &&& \
}
\]
The composition law \eqref{comp} is easily checked to be associative. We thus
obtain a category $\NC$ whose objects are vector spaces and whose morphisms
are NC power series. The identity morphism corresponding to a vector space $V$
is the power series $\id_V\colon V\to V$ given by $\id_1=\id_V$ and $\id_n=0$
for $n>1$.

\begin{lemma}
A NC power series $\phi\colon U\to V$ is an isomorphism iff the linear map
$\phi_1\colon U\to V$ is an isomorphism.
\end{lemma}

\begin{pf}
If $\phi_1$ is an isomorphism one can inductively solve the equation $\psi\circ\phi
=\id$ for $\psi_n\colon V^{\tensor n}\to V$. Similarly one can find $\psi'$ with $\phi
\circ \psi'=\id$. By general nonsense $\psi=\psi'$ is then an inverse for $\phi$. The
converse is obvious.
\end{pf}

\subsection{} 

Suppose $\phi\colon V\to V$ is a NC power series with $\phi_1=\id_V$. For any 
rooted plane tree $T$ with $n$ leaves, we can form a linear map
\[\phi_T\colon V^{\tensor n}\to V\]
by thinking of the leaves of $T$ as inputs and using the vertices of $T$ to compose
the tensors $\phi_k$. For example either of the trees above corresponds to the map
\[\phi_T(v_1,\ldots,v_6)=\phi_6(\phi_2(v_1,v_2),v_3,\phi_3(v_4,v_5,v_6)).\]

We shall only consider trees all of whose vertices have valency $\geq 3$
since the assumption $\phi_1=\id_V$ implies that vertices of valency $2$
do not contributes anything new.

\begin{thm}
\label{pr:inverses}
Let $\phi$ be a NC power series satisfying $\phi_1=\id_V$. Then $\phi^{-1}_1=
\id_V$ and, for $n>1$
\[\phi^{-1}_n(v_1,\ldots,v_n)=\sum_T (-1)^{|V(T)|} \phi_T(v_1,\ldots,v_n)\]
where the sum is over rooted plane trees  with $n$ leaves all of whose vertices have
valency $\geq 3$.
\end{thm}

\begin{pf}
It is enough to check that if one defines $\phi^{-1}$ by the given formula the composite
$\phi\circ \phi^{-1}$ is the identity. Clearly $(\phi\circ \phi^{-1})_1=\id_V$. For $n>1$,
expanding the composite $(\phi\circ \phi^{-1})_n$ gives a finite sum of signed terms
of the form $\phi_T(v_1,\ldots,v_n)$ for trees $T$ with $n$ leaves and vertices of valency
$\geq 3$.  Each such tree appears twice: once for the term in \eqref{comp} where
$k=1$, and once for a term with $k$ equal to the valency of the root vertex of $T$.
These two terms appear with opposite signs and hence cancel.
\end{pf}

\subsection{} 

Similarly to \ref{ss:NC series}, one can define a category $\CC$ whose objects
are complex vector spaces and whose morphisms are {\it commutative power
series} $\phi:U\to V$, that is sequences of linear maps $\phi_n:S^n U\to V$,
$n\geq 1$ where $S^n U$ is the $n$th symmetric power of $U$. Composition
is defined on tensors of the form $u^{\otimes n}$, $u\in U$, by
\begin{multline}
(\psi\circ\phi)_n(u,\ldots,u)=\\
\sum_{k=1}^n\thickspace\thickspace
\sum_{{0=i_0<\cdots <i_k=n}}
\psi_k \big(\phi_{i_1-i_0}(u,\ldots,u),\cdots,
\phi_{i_k-i_{k-1}}(u,\ldots,u)\big)
\label{eq:comm comp}
\end{multline}
and then by polarisation
$$(\psi\circ\phi)_n(u_1,\ldots,u_n)=
\frac{1}{n!}
\left.\frac{\partial}{\partial t_1}\right|_{t_1=0}
\negthickspace\negthickspace\negthickspace\negthickspace\negthickspace \cdots 
\left.\frac{\partial}{\partial t_n}\right|_{t_n=0}\,
(\psi\circ\phi)_n((t_1u_1+\cdots+t_nu_n)^{\otimes n})$$

The symmetrisation homomorphism $\sigma:S^nU\hookrightarrow U^{\otimes n}$
given by
$$\sigma(u_1\otimes\cdots\otimes u_n)=\frac{1}{n!}\sum_{\tau\in\Sym_n}
u_{\tau(1)}\otimes\cdots\otimes u_{\tau(n)}$$
yields a surjective restriction map $\sigma^*:\Hom_\NC(U,V)\to\Hom_\CC(U,V)$. A
comparison of \eqref {comp} and \eqref{eq:comm comp}, and the fact that a linear
map $S^nU\to V$ is uniquely determined by its values on tensors $u^{\otimes n}$,
$u\in U$, readily shows that $\sigma^*$ gives rise to a functor $\sigma^*:\NC\to\CC$.
In particular, Theorem \ref{pr:inverses} can be used to invert commutative power
series.

\subsection{}\label{ss:subspace}

If $U\subseteq V$ is a subspace, we shall say that a NC power series $\phi:V\to V$
{\it preserves} $U$ if the restriction $\sigma^*\phi_n$ of each $\phi_n$ to $S^n U
\subset V^{\otimes n}$ maps into $U$.

\begin{lemma}\label{le:subspace}\hfill
\begin{enumerate}
\item If $\phi,\psi$ preserve $U$, so does $\psi\circ\phi$.
\item If $\phi$ is invertible and preserves $U$, then so does $\phi^{-1}$.
\end{enumerate}
\end{lemma}
\begin{pf}
Denote the inclusion $U\hookrightarrow V$ by $\imath$ and let $p:V\to U$ be a
projection. A NC power series $\Theta$ preserves $U$ if, and only if $(\imath p-\id)
\circ\sigma^*\Theta\circ\imath=0$. (i) now follows since
$$(\imath p-\id)\circ\sigma^*(\psi\circ\phi)\circ i=
\left((\imath p-\id)\circ\sigma^*\psi\right)\circ\imath p\circ\sigma^*\phi\circ i=
\left((\imath p-\id)\circ\sigma^*\psi\circ\imath\right)\circ p\circ\sigma^*\phi\circ i=0$$
(ii) follows similarly from
$$0=
(\imath p-\id)\circ\sigma^*\id\circ i=
(\imath p-\id)\circ\sigma^*\phi^{-1}\circ\sigma^*\phi\circ\imath=
\left((\imath p-\id)\circ\sigma^*\phi^{-1}\circ\imath\right)\circ
\left( p\circ\sigma^*\phi\circ\imath\right)$$
and the fact that the commutative power series $ p\circ\sigma^*\phi\circ\imath:U\to U$
is invertible.
\end{pf}

\subsection{}

We shall be particularly interested in NC power series of the special form appearing in
Section \ref{se:Stokes map}. These power series depend on systems of coefficients
which we axiomatise as follows.

\begin{defn}
By a {\it transform} $F$ we mean a sequence of functions
\[F_n\colon \C^n\to \C,\quad  n\geq 1.\]
Given transforms $F$ and $G$ the composite transform $G\circ F$ is defined by the
finite sum
\begin{multline}
\label{composite}
(G\circ F)_n (z_1,\ldots,z_n)
=
\sum_{k=1}^n\thickspace
\sum_{0=i_0<\cdots <i_k=n} \Bigg[\;
G_k\bigg(\sum_{i=i_0+1}^{i_1}z_i, \sum_{i=i_1+1}^{i_2} z_i, \cdots, \sum_{i=i_{k-1}+1}^{i_k} z_i\bigg)\cdot
\\
\prod_{j=1}^k F_{i_{j}-i_{j-1}}(z_{i_{j-1}+1},\ldots, z_{i_{j}}) \;\Bigg].
\end{multline}
\end{defn}
Once again this sum is best thought of as a sum over trees of height 2. For example
the term corresponding to the tree \eqref{eq:tree} is
\[G_3(z_1+z_2,z_3,z_4+z_5+z_6) F_2(z_1,z_2) F_1(z_3) F_3(z_4,z_5,z_6).\]
The formula \eqref{composite} defines an associative composition law on the class of
transforms. The transform $\id$  with $\id_1=1$ and $\id_n=0$ for $n>1$ is  a two-sided
identity. It is easy to see that a transform $F$ is invertible precisely if the function $F_1$
is nowhere vanishing. Indeed, as before, in that case one can solve the equations $G
\circ F=\id$ and $F\circ H=\id$ inductively.

\subsection{}\label{ss:NC transform}

Transforms give rise to NC power series as follows. Let
\[A=\bigoplus_{\lambda\in\Lambda} A_{\lambda}\]
be an associative algebra over $\C$ graded by a free abelian group $\Lambda$.
For each $\lambda\in\Lambda$ let $\pi_\lambda\colon A\to A_\lambda$ be the
corresponding projection map. Suppose that we are given a fixed homomorphism
of abelian groups $Z\colon \Lambda\to\C$. Given a transform $F$ the corresponding
NC power series $\phi(F)\colon A\to A$ is given by the sum
\begin{equation}\label{transform}
\phi(F)_n(a_1,\ldots,a_n)=
\sum_{\lambda_1,\cdots,\lambda_n\in \Lambda}
F_n(Z(\lambda_1),\ldots,Z(\lambda_n))\,
\pi_{\lambda_1}(a_1)* \cdots *\pi_{\lambda_n}(a_n).
\end{equation}
It is  easy to check that given transforms $F$ and $G$ one has
\[\phi(G\circ F)=\phi(G)\circ\phi(F),\]
and hence that the data $(A,\Lambda,Z)$ defines a functor from the category
of transforms to that of NC power series $A\to A$.

\subsection{}

Suppose that $F$ is a transform satisfying $F_1(z)=1$ for all $z$. The method of
Theorem \ref{pr:inverses} allows us to give an explicit formula for the inverse of $
F$. We first associate a function
\[F_T\colon \C^n\to \C\]
to  a  plane rooted tree with $n$ leaves $T$ in the following way. Identify the leaves of
$T$ with their canonical order with the set $1,\ldots,n$. For each edge $e\in E(T)$ let
$I(e)\subset \{1,\ldots,n\}$ be the set of vertices lying above
$e$ and define the partial sum $s_e\colon \C^n\to \C$
\[s_e= s_e(z_1,\ldots,z_n)=\sum_{i\in I(e)} z_i.\]
To each vertex $v\in V(T)$  associate a factor
\[F_v(z_1,\ldots,z_n)=F_{m}\big(s_{e_1},s_{e_2},\ldots,s_{e_{m}}\big),\]
where $m+1$ is the valency of $v$ and $e_0,e_1,\ldots, e_m$
are the incident edges with their clockwise ordering, with $e_0$ being the outward pointing
edge. Then define $F_T$ to be the  product over vertices
\begin{equation}
\label{stanage}
F_T(z_1,\ldots,z_n)=\prod_{v\in V(T)} F_v(z_1,\ldots,z_n).\end{equation}
For example, for the tree $T$ depicted above
\[F_T(z_1,\ldots,z_6)=F_2(z_1,z_2) F_3(z_4,z_5,z_6) F_3(z_1+z_2,z_3,z_4+z_5+z_6).\]
The same argument as for Theorem \ref{pr:inverses} gives

\begin{prop}
\label{field}
Suppose $F$ is a transform satisfying $F_1(z)=1$ for all $z$. Then $F
^{-1}_1(z)=1$ for all $z$ and, for $n>1$ 
$$F^{-1}_n(z_1,\ldots,z_n)=\sum_T (-1)^{|V(T)|} F_T(z_1,\ldots,z_n)$$
where the sum is over rooted plane trees  with $n$ leaves all of whose
vertices have valency $\geq 3$.
\end{prop}

\subsection{}

The transforms most relevant to us have the following additional property.

\begin{defn}
A {\it Lie transform} is a transform $F$ such that, for any $n\geq 1$, $z_1,
\ldots,z_n\in\IC$ and non--commuting variables $x_1,\ldots,x_n$, the finite
sum
$$\sum_{\sigma\in\Sym_n}
F_n(z_{\sigma(1)},\ldots,z_{\sigma(n)})
x_{\sigma(1)}\cdots x_{\sigma(n)}$$
is a Lie polynomial in $x_1,\ldots,x_n$.
\end{defn}

\begin{lemma}\label{le:Lie transforms}\hfill
\begin{enumerate}
\item If $F$ and $G$ are Lie transforms, then so is $G\circ F$.
\item If $F$ is an invertible Lie transform, then so is $F^{-1}$.
\end{enumerate}
\end{lemma}
\begin{pf}
We begin by giving an alternative characterisation of a Lie transform $H$
in terms of the NC power series $\phi(H)$ introduced in \ref{ss:NC transform}.

For any $n\geq 1$, let $\L_n$ be the free Lie algebra on generators $x_1,
\ldots,x_n$. Its enveloping algebra $A=U\L_n$ possesses a grading by $
\Lambda=\IZ^n$ given by $\deg(x_i)=e_i$, where $\{e_1,\ldots,e_n\}$ is the
standard basis of $\IZ^n$. A transform $H$ is a Lie transform if, and only
if $\phi(H):U\L_n\to U\L_n$ preserves $\L_n\subset U\L_n$ in the sense
of \S \ref{ss:subspace} for every $n\geq 1$ and homomorphism $Z:\IZ^n\to
\IC$.

Indeed, if $H$ is a Lie transform, then for any $n\geq 1$, $x\in\L_n$ and
$m\geq 1$,
$$\phi(H)_m(x^{\otimes m})=
\sum_{\lambda_1,\ldots,\lambda_m}
H_m(Z(\lambda_1),\ldots,Z(\lambda_m))
\pi_{\lambda_1}(x)\cdots\pi_{\lambda_m}(x)$$
is a Lie polynomial in the variables $\pi_\lambda(x)$ and therefores lies
in $\L_n$. Conversely, if $\phi(H)$ preserves $\L_n$, and $Z:\IZ^n\to\IC$
is defined by $Z(e_i)=z_i$, the component of $\phi(H)_n((x_1+\cdots+x_n)
^{\otimes n})$ of weight $e_1+\cdots+e_n$ is equal to
$$\sum_{\sigma\in\Sym_n}
H_n(z_{\sigma(1)},\ldots,z_{\sigma(n)})
x_{\sigma(1)}\cdots x_{\sigma(n)}$$

The statements (i) and (ii) now follow from the foregoing and Lemma \ref
{le:subspace} since $\phi(G\circ F)=\phi(G)\circ\phi(F)$ and $\phi(F^{-1})
=\phi(F)^{-1}$.
\end{pf}

\subsection{}\label{ss:zero sum} 

We shall in fact need to consider transforms whose associated functions
$F_n:\IC^n\to\IC$ with $n\geq 2$ are only defined on $(\IC^*)^n$ and
satisfy
\begin{equation}\label{eq:zero sum}
F_n(z_1,\ldots,z_n)=0\qquad\text{whenever}\qquad
z_1+\cdots+z_n=0
\end{equation}
We will in such cases
tacitly extend the functions $F_n$ to $\IC^n$ in an arbitrary way. This does not
affect the values of the composition $G\circ F$ of two such transforms on $\Cstar$
since, by \eqref{composite} any evaluation of $G_k$ at an argument $w_{i_k}=
\sum_{i=i_{k-1}+1}^{i_k}z_i$ is multiplied by $F_{i_k-i_{k-1}}(z_{i_{k-1}+1},\ldots,
z_{i_k})$ which vanishes if $w_{i_k}$ does. Moreover, the composition $G\circ F$
satisfies \eqref{eq:zero sum} if $F$ and $G$ do.

Similarly, if $F$ is an invertible transform satisfying these properties, the inverse
transform $F^{-1}$ satisfies \eqref{eq:zero sum}, and the values of $(F^{-1})_n$
on $(\IC^*)^n$ do not depend on the extension of the functions $F_m$ to
$\IC^m$.
\smallskip

Let now $(A,\Lambda,Z)$ be a graded algebra as in \ref{ss:NC transform},
and consider the subspace $A^\times_Z\subset A$ defined by
$$A^\times_Z=
\bigoplus_{\substack{\lambda\in\Lambda\\Z(\lambda)\neq 0}}A_\lambda.$$
If $F$ is a transform defined on $\Cstar$, the restriction of $\phi(F)_n$ to $
(A^ \times_Z)^{\otimes n}$ is clearly independent of the extension of $F_n$
to $\IC^n$. Moreover, if $F$ satisfies \eqref{eq:zero sum}, the image of $\phi
(F)_n$ lies in $A^\times_Z$ for $n\geq 2$, so that $\phi(F)$ restricts to a NC
power series $A^\times_Z\to A^\times_Z$.

\subsection{} 

Let $U$ and $V$ be \fd and $\phi\colon U\to V$ a commutative power series.
Assume that the sum
\begin{equation}
\label{taylor}
\underline{\phi}(u)=\sum_{n\geq 1}\phi_n(u^{\tensor n})
\end{equation}
is  convergent for all $u$ in an open neighbourhood of the origin $0\in U^o
\subset U$. Then  $\underline{\phi}$ defines a holomorphic map $U^o\to V$
and \eqref{taylor} is its Taylor expansion at the origin.

If $\phi:U\to V$ and $\psi\colon V\to W$ are two commutative power series
which convergent in neighbourhoods of the origins in $U$ and $V$ respectively
then the power series $\psi\circ\phi$ is convergent in a neighbourhood of the
origin in $U$ and $\underline{(\psi\circ\phi)}=\underline{\psi}\circ\underline{\phi}$.

Similarly, if $V$ is \fd and $\phi\colon V\to V$ is an invertible power series
which is convergent in a neighbourhood of the origin, then by the inverse
function theorem, the inverse map $\underline{\phi}^{-1}$ is holomorphic
near the origin.  The following is standard.

\begin{lemma}\label{le:sum of inverse}
The power series $\phi^{-1}$ is convergent in a \nbd of the origin in $V$
and $\ul{\phi^{-1}}={\ul{\phi}}^{-1}$.
\end{lemma}
\begin{pf}
If $f:V\to V$ is a germ of a holomorphic function at $0\in V$ such that $f(0)=0$,
we denote its Taylor series, viewed as a commutative power series, by $Tf$.
Thus, $\ul{Tf}=f$ and $T(\ul{\phi})=\phi$ whenever $\ul{\phi}$ is defined. Since
$\ul{\phi}\circ\ul{\phi}^{-1}=\id_V=\ul{\phi}^{-1}\circ\ul{\phi}$ we find, upon applying
$T$ that 
$$\phi\circ T(\ul{\phi}^{-1})=\id_V=T(\ul{\phi}^{-1})\circ\phi$$
Thus $\phi^{-1}=T(\ul{\phi}^{-1})$ as claimed.
\end{pf}

\subsection{Proof of Theorem \ref{th:Stokes inverse}}

We shall proceed by extending the Taylor series of the Stokes map
$\calS$ to a commutative power series $U\g\to U\g$, then lift it to a
NC Lie transform $\phi(L):U\g\to U\g$ and finally invert it by using
Theorem \ref{pr:inverses}. The Taylor series of $\calS^{-1}$ will
then be obtained as the restriction to $\god\subset U\g$ of the
commutative power series $\sigma^*\phi(L)^{-1}$.

Specifically, let $U\g$ be the universal enveloping algebra of $\g$, graded
by the lattice $\Lambda\subset\h^*$ spanned by the set of roots $\Phi(G;H)$.
Let $L_n:(\IC^*)^n\to\IC$ be the functions defined in \ref{ss:Stokes map},
$L=\{L_n\}_{n\geq 1}$ the corresponding transform and $\phi(L):U\g\to
U\g$ the NC power series determined by $(U\g,\Lambda,Z)$ and $L$.
Thus, for any $n\geq 1$ and $x^1,\ldots,x^n\in U\g$,
\begin{equation}\label{eq:lift}
\phi(L)_n(x^1\otimes\cdots\otimes x^n)=
\sum_{\gamma_1,\ldots,\gamma_n}
L_n(Z(\gamma_1),\ldots,Z(\gamma_n))\medspace
x^1_{\gamma_1}\cdots x^n_{\gamma_n}
\end{equation}
where $x_\gamma$ is the weight component of $x$ corresponding to
$\gamma\in\Lambda$.

The power series $\phi(L)$ preserves the subspace $\g\subset U\g$ in the
sense of \ref{ss:subspace} since $L$ is a Lie transform by Theorem \ref{one}.
By \ref{ss:zero sum}, each $\phi(L)_n$, $n\geq 2$, maps $(U\g)^{\otimes n}$
to
$$U\g^\times=
\bigoplus_{\substack{\lambda\in\Lambda\\ Z(\lambda)\neq 0}}(U\g)_\lambda$$
since $L$ satisfies \eqref{eq:zero sum} by \eqref{eq:L zero sum}. Thus,
$\phi(L)_n$ maps $S^n\god$ to $\god=U\g^\times\cap\g$ so that the
restriction of $\sigma^*\phi(L)$ to $\god\subset U\g$ is a commutative
power series $\god\to\god$ which is equal to the Taylor series of the
Stokes map $\calS$ by Theorem \ref{one}.

Let now $J=L^{-1}$ be the inverse transform. By Proposition \ref{field},
$J_1\equiv (2\pi i)^{-1}$ and, for $n\geq 2$,
$$J_n(z_1,\ldots,z_n)=(2\pi i)^{-n}\sum_T (-1)^{|V(T)|} J_T(z_1,\ldots,z_n)$$
where the sum is over rooted plane trees  with $n$ leaves all of whose
vertices have valency $\geq 3$ and $J_T$ is defined by \eqref{stanage}.
$J$ is a Lie transform by Lemma \ref{le:Lie transforms} which satisfies
\eqref{eq:zero sum} since $L$ does. Thus $\sigma^*\phi(J)=\sigma^*\phi
(L)^{-1}$ restricts to a commutative power series $\god\to\god$. The
latter is given by \eqref{main2} and is equal to the Taylor series of $\calS
^{-1}$ at $\epsilon=0$ by Lemma \ref{le:sum of inverse}.
\qed

\remark The transform $\phi(L)$ is in a sense the most economical lift of
the Taylor series of $\calS$ to a NC power series. In particular, it differs
from the more canonical lift obtained by composing the terms $S^n\god
\to\god$ of the Taylor series of $\calS$ with the canonical projection $\god
^{\otimes n}\to S^n\god$.



\begin{thebibliography}{101}

\bibitem{BJL} W. Balser, W. B. Jurkat, and D. A. Lutz, {\it Birkhoff invariants
and StokesÕ multipliers for meromorphic linear differential equations}, J.
Math. Anal. Appl. {\bf 71} (1979), 48--94.

\bibitem{BJL2} W. Balser, W. B. Jurkat, and D. A. Lutz, {\it On the reduction
of connection problems for differential equations with an irregular singular
point to ones with only regular singularities, I}, SIAM J. Math. Anal. {\bf 12}
(1981), 691--721.

\bibitem{BJL3} W. Balser, W. B. Jurkat, and D. A. Lutz, {\it Characterization
of first level formal solutions by means of the growth of their coefficients},
J. Differential Equations {\bf 51} (1984), no. 1, 48--77. 
 
\bibitem{BJL4} W. Balser, W. B. Jurkat, and D. A. Lutz, {\it Transfer of
connection problems for first level solutions of meromorphic differential
equations, and associated Laplace transforms}, J. Reine Angew. Math.
{\bf 344} (1983), 149--170.

\bibitem{Boalch} P. P. Boalch, {\it Symplectic manifolds and isomonodromic
deformations}, Adv. Math. {\bf 163} (2001), 137--205.

\bibitem{Boalch3} P. P. Boalch, {\it Stokes matrices, Poisson Lie groups
and Frobenius manifolds}, Invent. Math. {\bf 146} (2001), no. 3, 479--506.

\bibitem{Boalch2} P. P. Boalch, {\it $G$--bundles, isomonodromy, and
quantum Weyl groups}, Int. Math. Res. Not. {\bf 2002}, 1129--1166.

\bibitem{BTL} T. Bridgeland and V. Toledano Laredo, {\it Stability
conditions and Stokes factors}, {\sf arXiv:0801.3974}, to appear
in Invent. Math. (2011).

\bibitem{Gon} A. B. Goncharov, {\it Multiple polylogarithms and mixed
Tate motives}, {\sf math.AG/0103059}.

\bibitem{Hain} R. M. Hain, {\it The geometry of the mixed Hodge structure
on the fundamental group}. Algebraic geometry, Bowdoin, 1985 (Brunswick,
Maine, 1985), 247--282, Proc. Sympos. Pure Math., 46, Part 2, AMS, 1987.

\bibitem{Hardy} G. H. Hardy, {\it Divergent series}. Clarendon Press, Oxford,
1949.

\bibitem{Hu} J. E. Humphreys, {\it Linear algebraic groups}. Graduate Texts
in Mathematics, No. 21. Springer--Verlag, 1975.

\bibitem{Ince} E. L. Ince, {\it Ordinary Differential Equations}. Dover Publications,
1944.

\bibitem{JMU} M. Jimbo, T. Miwa and K. Ueno, {\it Monodromy preserving
deformation of linear ordinary differential equations with rational coefficients.
I.} Phys. D {\bf 2} (1981), 306--352.

\bibitem{Joyce} D. Joyce, {\it Holomorphic generating functions for invariants
counting coherent sheaves on Calabi--Yau 3--folds}, Geom. Topol. {\bf 11}
(2007), 667--725.

\bibitem{JLP} W. Jurkat, D. Lutz, and A. Peyerimhoff, {\it Birkhoff invariants
and effective calculations for meromorphic linear differential equations 1.},
J. Math. An. Appl. {\bf 53} (1976), 438--470.

\bibitem{Reineke}  M. Reineke, {\it The Harder--Narasimhan system in
quantum groups and cohomology of quiver moduli}, Invent. Math. {\bf
152} (2003), 349--368.

\bibitem{Rudin} W. Rudin, {\it Function theory in the unit ball of $\C\sp{n}$}.
Grundlehren der Mathematischen Wissenschaften, 241. Springer--Verlag,
1980.

\bibitem{SR} N. Saavedra Rivano, {\it Cat\'egories Tannakiennes}.
Lecture Notes in Mathematics, Vol. 265. Springer--Verlag, 1972.

\bibitem{Wasow} W. Wasow, {\it Asymptotic expansions for ordinary
differential equations}. Dover Publications, 1987.

\end{thebibliography}
\end{document}